\documentclass[MSNbibl,number,citesort,seceqn,dvips]{arxbj}
\usepackage{upgreek}
\usepackage{graphicx}

\aid{0}
\volume{19}
\issue{5A}
\pubyear{2013}
\firstpage{1655}
\lastpage{1687}
\doi{10.3150/12-BEJ425} 

\makeatletter
\newtheorem{lemma}{Lemma}[section]
\newtheorem{theorem}[lemma]{Theorem}
\newremark{rem}[lemma]{Remark}
\newtheorem{proposition}[lemma]{Proposition}

\newcommand{\eqref}[1]{(\ref{#1})}
\newcommand{\Nat}{\mathbb{N}}
\newcommand{\Q}{\mathbb{Q}}
\newcommand{\R}{\mathbb{R}}
\newcommand{\Z}{\mathbb{Z}}
\newcommand{\eps}{\varepsilon}
\newcommand{\vect}[1]{\mathbf{#1}} 
\newcommand{\Prob}{\mathbb{P}} 
\newcommand{\Exp}{\mathbb{E}} 
\newcommand{\ind}{\mathbb{I}} 
\newcommand{\Ec}{\mathcal{E}}
\newcommand{\Bc}{\mathcal{B}}
\newcommand{\Db}{\mathbb{D}}
\newcommand{\Gb}{\mathbb{G}}
\newcommand{\thetamd}{\hat\theta_n^{\mathrm{MD}}}
\newcommand{\id}{\mathrm{id}}
\newcommand{\ran}{\operatorname{ran}}

\newcommand{\ltc}{ \Lambda_{L} }
\newcommand{\ltcein}{ \Lambda_{L}^\angle}
\newcommand{\ltchat}{ \hat{ \Lambda}_{L} }
\newcommand{\ltchatein}{ \hat{ \Lambda}_{L}^\angle}
\newcommand{\ltctilde}{\tilde{\Lambda}_{L}}
\newcommand{\ltci}[1]{ \Lambda_{L, #1 } }
\newcommand{\ltciein}[1]{ \Lambda_{L, #1 }^\angle}
\newcommand{\ltcihat}[1]{ \hat{ \Lambda}_{L, #1 } }
\newcommand{\ltcihatein}[1]{ \hat{ \Lambda}_{L, #1 }^\angle}

\newcommand{\utc}{ \Lambda_{U} }
\newcommand{\utchat}{ \hat{ \Lambda}_{U} }
\newcommand{\utctilde}{\tilde{\Lambda}_{U}}

\newcommand\weak{\rightsquigarrow}
\newcommand{\weakP}{ \mathop{\stackrel{\Prob}{\rightsquigarrow}}_{\xi} }
\newcommand{\Pconv}{ \stackrel{\Prob}{\rightarrow} }
\makeatother

\begin{document}
\begin{frontmatter}

\title{Multiplier bootstrap of tail copulas with applications}
\runtitle{Multiplier bootstrap of tail copulas}

\begin{aug}
\author{\fnms{Axel} \snm{B\"{u}cher}\thanksref{e1}\ead[label=e1,mark]{axel.buecher@ruhr-uni-bochum.de}} \and
\author{\fnms{Holger} \snm{Dette}\corref{}\thanksref{e2}\ead[label=e2,mark]{holger.dette@ruhr-uni-bochum.de}}
\runauthor{A. B\"{u}cher and H. Dette} 
\address{Fakult\"at f\"ur Mathematik, Ruhr-Universit\"{a}t Bochum,
44780 Bochum, Germany.\\ \printead{e1,e2}}
\end{aug}

\received{\smonth{1} \syear{2011}}
\revised{\smonth{1} \syear{2012}}

%
\begin{abstract}
For the problem of estimating lower tail and upper tail copulas, we propose
two bootstrap procedures for approximating the distribution of the corresponding
empirical tail copulas. The first method uses a multiplier bootstrap of
the empirical
tail copula process and requires estimation of the partial derivatives
of the tail
copula. The second method avoids this estimation problem and uses multipliers
in the two-dimensional empirical distribution function and in the
estimates of the
marginal distributions. For both multiplier bootstrap procedures, we prove
consistency.

For these investigations, we demonstrate that the common assumption of
the existence of continuous partial derivatives in the
the literature on tail copula estimation is so restrictive, such that
the tail copula corresponding to tail independence is the only tail
copula with this property.
Moreover, we are able to solve this problem and prove weak convergence
of the empirical tail copula process under nonrestrictive
smoothness assumptions that are satisfied for many commonly used
models. These
results are applied in several statistical problems, including
minimum distance estimation and goodness-of-fit testing.
\end{abstract}

%

\end{frontmatter}

\section{Introduction}
The stable tail dependence function appears naturally in multivariate
extreme value theory as a function that characterizes extremal
dependence. If a bivariate distribution function $F$ lies in the
max-domain of attraction
of an extreme-value distribution $G$, then the copula of $G$ is completely
determined by the stable tail dependence function (see, e.g., Einmahl \textit{et~al.}
\cite{einkraseg2008}).
The function is closely related to tail copulas (see, e.g., Schmidt and Stadtm{\"{u}}ller \cite{schmstad2006})
and represents the current standard to describe extremal
dependence (see Embrechts \textit{et~al.} \cite{emblinmcn2003} and Malevergne and Sornette \cite{malsor2004}).
The lower and upper tail copulas are defined by
%
\begin{equation} \label{lp}
\ltc(\vect{x}) = \lim_{t\rightarrow\infty} t C(x_1/t,x_2/t),\qquad
\utc(\vect{x}) = \lim_{t\rightarrow\infty} t \bar{C}(x_1/t,x_2/t),
\end{equation}
provided that the limits exist. Here
$\vect{x}=(x_1,x_2)\in\bar{\R}_+^2:=[0,\infty]^2\setminus\{
(\infty,\infty)\}$,
$C$ denotes the (unique) copula of the two-dimensional continuous
distribution function $F$,
which relates $F$ and its marginals $F_1, F_2 $ by
%
\begin{equation}
F(\vect{x}) = C(F_1(x_1),F_2(x_2)) \label{cop}
\end{equation}
(see Sklar \cite{sklar1959}), and
$
\bar C(\vect{u})=u_1+u_2-1+C(1-u_1,1-u_2) 
$
denotes the survival copula of $\vect{X}=(X_1,X_2)\sim F$. The stable
tail dependence function $l$ and the
upper tail copula $\utc$ are associated through the relationship
\[
l(\vect{x})=x_1+x_2 - \utc(\vect{x}) \qquad \forall\vect{x}\in\R_+^2.
\]
Since its introduction various parametric and nonparametric estimates
of the tail copulas and of the stable tail dependence function have
been proposed in the literature. Several authors assume that the
dependence function belongs to some parametric family. Coles and Tawn \cite{coletawn1994},
Tiago~de Oliveira \cite{tiago1980}, and Einmahl \textit{et~al.} \cite{einhaaxin1993} imposed restrictions on
the marginal distributions to estimate multivariate extreme value distributions.
Nonparametric estimates of the stable tail dependence function were investigated
in the pioneering thesis of Huang \cite{huang1992} and by Qi \cite{qi1997}
and Drees and Huang \cite{dreehuan1998}. Schmidt and Stadtm{\"{u}}ller \cite{schmstad2006} proposed analogous
estimates as in Huang \cite{huang1992} for tail copulas [except for
rounding deviations due to the fact that $F_n(F_n^-(x))$, with the
generalized inverse function $F_n^-$, is not exactly equal to $x$] and
provided new proof of the asymptotic behavior of the estimates. More
recent work on inference on the stable tail dependence function was
done by Einmahl \textit{et~al.} \cite{einkraseg2008} and Einmahl \textit{et~al.} \cite{einhaali2006}, who
investigated moment estimators of tail dependence and weighted
approximations of tail copula processes, respectively.

The present paper has two main purposes. First, we clarify some
curiosities in the literature
on tail copula estimation, which stem from the fact that most authors
assume the existence of continuous
partial derivatives of the tail copula {(see, e.g., Huang \cite{huang1992},
Drees and Huang \cite{dreehuan1998}, Schmidt and Stadtm{\"{u}}ller \cite{schmstad2006}, Einmahl \textit{et~al.} \cite{einhaali2006},
de~Haan and Ferreira \cite{dehaferr2006}, Peng and Qi \cite{pengqi2008}, de~Haan \textit{et~al.} \cite{dehnevpen2008}).}
{However, the (lower or upper) tail copula corresponding to (lower or
upper) tail independence is the only tail copula with this property,
because the
partial derivatives of a tail copula satisfy
%
\begin{equation} \label{counter}
\partial_1 \Lambda(0,x) =
\cases{
\displaystyle\lim_{t\rightarrow\infty} \Lambda(1,t) &\quad if $x > 0$,\vspace*{2pt}
\cr
0 &\quad if $x =0$,}
\end{equation}
where $\Lambda$ denotes either $\ltc$ or $\utc$ (see Appendix~\ref
{app:B} for details).}
Consequently, we provide a result on the weak convergence of the
empirical tail copula process
(and thus also of the
empirical stable tail dependence function) under weak smoothness
assumptions (see Theorem~\ref{theo:ltchat} in the next section). The
smoothness conditions are nonrestrictive in the sense that in the case where
they are not satisfied, the candidate limiting process does not have
continuous trajectories.

Note that similar investigations were recently carried out by
Segers \cite{segers2011} in the context of nonparametric copula estimation. In that
paper it is demonstrated that many (even most) of the most popular
copula models do not have continuous partial derivatives on the whole
unit square, which has been the usual assumption for the asymptotic
behavior of the empirical copula process hitherto. Moreover, it is
shown how the assumptions can be suitably relaxed such that the
asymptotics are not influenced.

The second objective of the present paper is to approximate the distribution
of estimators for the tail copulas by new bootstrap methods. In
contrast to the
problem of estimation of the stable dependence function and tail copulas,
this problem has received much less attention in the literature. Recently,
Peng and Qi \cite{pengqi2008} considered the tail empirical distribution function and
showed the consistency of the bootstrap based on resampling (again
under the assumption of continuous partial
derivatives). They used their results to construct confidence bands for
the tail dependence
function. Although the authors considered the naive bootstrap, the present
paper is devoted to multiplier bootstrap procedures for tail copula estimation.
On the one hand, our research is motivated by the observation that the
parametric bootstrap, which is commonly applied in goodness-of-fit testing
problems (see de~Haan \textit{et~al.} \cite{dehnevpen2008}), has very high computational costs
because it relies heavily on random number generation and estimation.
(See Kojadinovic and Yan \cite{kojayan2010} and Kojadinovic \textit{et~al.} \cite{kojyanhol2010} for a more
detailed discussion of the
computational efficiency of the multiplier bootstrap.)
On the other hand, as pointed out by B{\"{u}}cher \cite{buecher2011} and
B{\"{u}}cher and Dette \cite{buecdett2010} in the context
of nonparametric copula estimation, some multiplier bootstrap procedures
lead to more reliable approximations compared with those from the
bootstrap based on resampling.

In Section~\ref{sec:etc} we briefly review the nonparametric estimates
of the
tail copula and discuss their main properties. In particular, we establish
weak convergence of the empirical tail copula process under
nonrestrictive smoothness
assumptions, which are satisfied for many commonly used models.
In Section~\ref{sec:tcmult} we introduce the multiplier bootstrap for
the empirical
tail copula and prove its consistency. In particular, we discuss two
ways of
approximating the distribution of the empirical tail copula by a
multiplier bootstrap.
Our first method, called the \textit{partial derivatives multiplier
bootstrap}, uses the structure of the
limit distribution of the empirical tail copula process. As a
consequence, this
approach requires estimation of the partial derivatives of the tail copula.
The second method, which we call the \textit{direct multiplier
bootstrap}, avoids this problem, using multipliers in the
two-dimensional empirical distribution
function and in the estimates of the marginal distributions.
Finally, in Section
\ref{sec:tcapplications} we discuss several statistical applications
of the
multiplier bootstrap. In particular, we investigate the problem of
testing for
equality between two tail copulas and discuss the bootstrap approximations
in the context of testing parametric assumptions for the tail copula.
We defer all proofs and some of the technical details to the \hyperref[append]{Appendix}.

\section{Empirical tail copulas}\label{sec:etc}
Let $\vect{X}_1, \ldots, \vect{X}_n$ denote independent identically
distributed random variables
with distribution function $F$ and denote the empirical distribution
functions of
$F$ and its marginals $F_1$ and $F_2$ by $F_n(\vect{x})=n^{-1}\sum
_{i=1}^n \ind\{ \vect{X}_i\leq\vect{x} \}$,
$F_{n1}(x_1)=F_{n}(x_1,\infty)$ and $F_{n2}(\vect{x})=F_n(\infty
,x_2)$, respectively. Analogously, we define the joint empirical
survival function by $\bar F_n(\vect{x}) = n^{-1}\sum_{i=1}^n \ind\{
\vect{X}_i > \vect{x} \}$ and the marginal empirical survival
functions as $\bar F_{n1}= 1- F_{n1}$ and $\bar F_{n2}=1-F_{n2}$. Following
Schmidt and Stadtm{\"{u}}ller \cite{schmstad2006}, we consider the estimators
%
\begin{equation}
\label{empl}
\ltchat(\vect{x}) = \frac{n}{k}C_n \biggl(\frac{kx_1}{n},\frac{kx_2}{n} \biggr),\qquad
\utchat(\vect{x}) = \frac{n}{k}\bar C_n\biggl (\frac{kx_1}{n},\frac
{kx_2}{n} \biggr)
\end{equation}
for the lower and upper tail copulas, respectively, where $k\rightarrow
\infty$ such
that $k=\mathrm{o}(n)$ and $C_n$ (resp., $\bar{C}_n$) denotes the empirical
copula (resp.,
empirical survival copula), that is,
%
\[
C_n(\vect{u}) =F_n(F_{n1}^-(u_1),F_{n2}^-(u_2)),\qquad
\bar C_n(\vect{u}) = \bar F_n(\bar F_{n1}^-(u_1),\bar F_{n2}^-(u_2)).
\]
Here $G^-$ and $\bar G^-$ denote the (left-continuous) generalized
inverse functions of some real distribution function $G$ and its
corresponding survival function $\bar G=1-G$, defined by
\begin{eqnarray*} G^{-}(p)&:=&
\cases{
\inf\{x\in\R\vert G(x)\geq p \}, &\quad $0<p\leq1$,
\cr
\sup\{x\in\R\vert G(x)=0\}, & \quad$p=0$,}
\\
\bar G^{-}(p)&:=&
\cases{
\sup\{x\in\R\vert \bar G(x)\geq p \}, & \quad$0<p\leq1$,
\cr
\inf\{x\in\R\vert \bar G(x)=0\}, & \quad$p=0$.}
\end{eqnarray*}

It is easy to see that the estimators $\ltchat$ and $\utchat$ are
asymptotically equivalent to the estimates
%
\begin{eqnarray}
\label{est:ltchat}\frac{1}{k}\sum_{i=1}^{n} \ind\{R(X_{i1})\leq kx_1, R(X_{i2})\leq
kx_2\}&=& \ltchat(\vect{x}) +\mathrm{O}(1/k) ,
\\
\label{est:utchat}\frac{1}{k}\sum_{i=1}^{n} \ind\{R(X_{i1})>n-kx_1,
R(X_{i2})>n-kx_2\} &=& \utchat(\vect{x}) +\mathrm{O}(1/k) ,
\end{eqnarray}
where $R(X_{ij})=nF_{n1}(X_{j1})$ denotes the rank of $X_{ij}$ among
$X_{1j},\ldots, X_{nj}$ ($j=1,2$) (see Huang \cite{huang1992} for an
alternative asymptotic equivalent estimator).
Therefore, we introduce analogs of (\ref{est:ltchat}) and (\ref
{est:utchat}) where
the marginals $F_{1}$ and $F_{2}$ are assumed known, that is,
%
\begin{eqnarray}
\label{est:ltctilde}
\ltctilde(\vect{x}) & =& \frac{1}{k}\sum_{i=1}^{n} \ind\biggl\{
F_1(X_{i1})\leq\frac{kx_1}{n}, F_2(X_{i2})\leq\frac{kx_2}{n}\biggr\},
\\
\label{est:utctilde}\utctilde(\vect{x}) &=& \frac{1}{k}\sum_{i=1}^{n} \ind\biggl\{
F_1(X_{i1})>1-\frac{kx_1}{n}, F_2(X_{i2})>1-\frac{kx_2}{n}\biggr\}.
\end{eqnarray}
For the sake of brevity, we restrict our investigations to the case of
lower tail
copulas. We assume that this function is non-zero in a single point
$\vect{x}\in[0,\infty)^2\setminus\{(0,0)\}$, and as a consequence
non-zero everywhere on
$[0,\infty)^2$ (see Theorem 1 in Schmidt and Stadtm{\"{u}}ller~\cite{schmstad2006}).

Let $\mathcal{B}_\infty(\bar{\R}_+^2)$ denote the space of all functions
$f\dvtx \bar{\R}_+^2\rightarrow\R$, which are locally uniformly bounded
on every
compact subset of $\bar{\R}^2_+=[0,\infty]^2\setminus\{(0,0)\}$
(i.e., on closed subsets\vadjust{\goodbreak} of $[0,\infty]^2$ that are bounded away from $0$),
equipped with the metric
\[
d(f_1,f_2) = \sum_{i=1}^\infty2^{-i} (\Vert f_1-f_2\Vert_{T_i}\wedge1),
\]
where the sets $T_i$ are defined by $T_i=[0,i]^2\cup[0,i]\times\{
\infty\}\cup\{\infty\}\times[0,i]$
and where\vspace*{1pt} $\Vert f\Vert_{T_i}=\sup_{\vect{x}\in T_i} |f(\vect{x})|$ denotes
the sup-norm on $T_i$. Note that with this metric, the set $\mathcal
{B}_\infty(\bar{\R}_+^2)$
is a complete metric space and that a sequence $f_n$ in $\mathcal
{B}_\infty(\bar{\R}_+^2)$
converges with respect to $d$ if and only if it converges uniformly on
every $T_i$
(see {V}an~der Vaart and Wellner \cite{vandwell1996}). Throughout\vspace*{-1pt} this paper, $\ell^{\infty}(T)$
denotes the set of
uniformly bounded functions on a set $T$, $\Pconv$ denotes convergence
in (outer) probability, and $ \weak$ denotes weak
convergence in the sense of Hoffmann-J\o rgensen (see, e.g.,
{V}an~der Vaart and Wellner \cite{vandwell1996}).

Schmidt and Stadtm{\"{u}}ller \cite{schmstad2006} assumed that the lower tail copula $\ltc$
satisfies the
second-order condition
%
\begin{equation} \label{secor}
\lim_{t\rightarrow\infty}\frac{\ltc(\vect
{x})-tC(x_1/t,x_2/t)}{A(t)}=g(\vect{x})
\end{equation}
locally uniformly for $\vect{x}=(x_1,x_2)\in\bar{\R}_+^2$, where
$g$ is a
non-constant function and the function $A\dvtx[0,\infty)\rightarrow
[0,\infty)$ satisfies
$\lim_{t\rightarrow\infty}A(t)=0$. A detailed look at the proofs
reveals that \eqref{secor} may be weakened to the condition
%
\begin{equation}\label{secornew}
| \ltc(\vect{x})-tC(x_1/t,x_2/t)| =\mathrm{O}(A(t))
\end{equation}
for $t\to\infty$, locally uniformly for $\vect{x}\in\bar{\R
}_+^2$. Under \eqref{secor} and the additional assumptions
$\ltc\not\equiv0$, $\sqrt{k}A(n/k)\rightarrow0 $,
$k=k(n)\rightarrow\infty$, $k=\mathrm{o}(n)$,
Schmidt and Stadtm{\"{u}}ller \cite{schmstad2006} showed that the lower tail copula process with
known marginals defined by (\ref{est:ltctilde}) converges
weakly in $\mathcal{B}_\infty(\bar{\R}_+^2)$, that is,
%
\begin{equation} \label{t4ss}
\sqrt{k} \bigl( \ltctilde(\vect{x})-\ltc(\vect{x}) \bigr) \weak\Gb
_{\ltctilde}(\vect{x}) ,
\end{equation}
where $\Gb_{\ltctilde}$ is a centered Gaussian field with covariance structure
given by
%
\begin{equation} \label{covgtilde}
\Exp \Gb_{\ltctilde} (\vect{x}) \Gb_{\ltctilde}(\vect{y}) = \ltc
(x_1\wedge y_1, x_2 \wedge y_2).
\end{equation}
For the empirical tail copula $\ltchat(\vect{x})$ the authors
established the weak convergence
%
\begin{equation} \label{t5ss}
\alpha_n(\vect{x}) = \sqrt{k}\bigl ( \ltchat(\vect{x})-\ltc(\vect{x})
\bigr) \weak\Gb_{\ltchat}(\vect{x})
\end{equation}
in $\mathcal{B}_\infty(\bar{\R}_+^2)$, provided that the tail
copula has
continuous partial derivatives. Here the limiting process $\Gb
_{\ltchat}$
has the representation
%
\begin{equation} \label{gf}
\Gb_{\ltchat}(\vect{x})= \Gb_{\ltctilde}(\vect{x})-\partial_1
\ltc(\vect{x}) \Gb_{\ltctilde}(x_1,\infty) - \partial_2 \ltc
(\vect{x}) \Gb_{\ltctilde}(\infty,x_2).
\end{equation}

The assumption of continuous partial derivatives is made in the
literature on estimation of stable tail dependence functions
and tail copulas. However, as demonstrated in \eqref{counter}, there
does not exist any
tail copula $\ltc\not\equiv0$ with continuous partial derivatives at
the origin $(0,0)$.
With our first result, we fill this gap and prove weak convergence of
the empirical tail copula process under
suitable weakened smoothness assumptions. To do so, we use a similar
approach as in Schmidt and Stadtm{\"{u}}ller \cite{schmstad2006}, because this
turns out to be useful for a proof of consistency of the multiplier
bootstrap as well.
We first consider the case of known marginals. Because of the
second-order condition (\ref{secornew}), the proof of (\ref{t4ss})
can be given by
showing weak convergence of the centered statistic
%
\begin{equation} \label{eq:atilde}
\tilde\alpha_n(\vect{x}) := \sqrt{k} \biggl(\ltctilde(\vect{x})-\frac
{n}{k} C(x_1k/n, x_2k/n) \biggr).
\end{equation}

\begin{lemma} \label{lem:ltctilde}
If $\ltc\not\equiv0$ and the second-order condition \eqref
{secornew} holds with
$\sqrt{k}A(n/k)\rightarrow0 $, where $k=k(n)\rightarrow\infty$ and
$k=\mathrm{o}(n)$, then we have, as $n$ tends to infinity,
%
\begin{equation}\label{eq:tildeakonv}
\tilde\alpha_n(\vect{x}) =\sqrt{k} \biggl(\ltctilde(\vect{x})-\frac
{n}{k} C(x_1k/n, x_2k/n) \biggr) \weak\Gb_{\ltctilde}(\vect{x})
\end{equation}
in $\mathcal{B}_\infty(\bar{\R}_+^2)$. Here $\Gb_{\ltctilde}$ is a
tight centered Gaussian field concentrated on $\mathcal{C}_\rho(\bar
{\R}^2_+)$
with covariance structure given in (\ref{covgtilde})
and $\rho$ is a pseudo-metric
on the space $\bar{\R}_+^2$ defined by
\[
\rho(\vect{x} ,\vect{y}) = \Exp\bigl[\bigl(\Gb_{\ltctilde}(\vect{x})-\Gb
_{\ltctilde}(\vect{y})\bigr)^2 \bigr]^{1/2}= \bigl( {\ltc}(\vect{x})-2{\ltc
}(\vect{x}\wedge\vect{y})+{\ltc}(\vect{x})
\bigr)^{1/2},
\]
$\vect{x}=(x_1,x_2), \vect{y}=(y_1,y_2), \vect{x}\wedge\vect{y} =
(x_1 \wedge y_1, x_2 \wedge y_2)$
and $\mathcal{C}_\rho( \bar{\R}^2_+ )\subset\mathcal{B}_\infty(
\bar{\R}^2_+ )$
denotes the subset of all functions that are uniformly $\rho
$-continuous on
every $T_i$.
\end{lemma}

This assertion is proved in Theorem 4 of Schmidt and Stadtm{\"{u}}ller \cite{schmstad2006} by
showing convergence of the finite-dimensional distributions and
tightness. For an alternative proof based on Donsker classes, see
Remark~\ref{rem:ap} in the \hyperref[appa]{Appendix}.
For a proof of a corresponding result for the empirical tail copula
process with estimated marginals as defined in \eqref{t5ss}, we use
the functional delta method in \eqref{t4ss} with some suitable functional.
\begin{theorem}\label{theo:ltchat}
Let $\ltc\not\equiv0$ be a lower tail copula whose first-order
partial derivatives satisfy the condition
%
\begin{equation}\label{fotc}
\partial_p \ltc\mbox{ exists and is continuous on } \{ \vect{x}\in
\bar\R_+^2 \vert 0 < x_p<\infty\}
\end{equation}
for $p=1,2$. If in addition the assumptions of Lemma~\ref
{lem:ltctilde} are satisfied, then we have
\[
\alpha_n(\vect{x}) = \sqrt{k} \bigl( \ltchat(\vect{x})-\ltc(\vect{x})
\bigr) \weak\Gb_{\ltchat}(\vect{x})
\]
in $\mathcal{B}_\infty(\bar{\R}_+^2)$, where the process $\Gb
_{\ltchat}$ is defined in \eqref{gf} and $\partial_p \ltc, p=1,2$
is defined as 0 on the set $\{\vect{x}\in\bar\R_+^2 \vert x_p\in\{
0,\infty\} \}$.
\end{theorem}

Theorem~\ref{theo:ltchat} has been proven by Schmidt and Stadtm{\"{u}}ller \cite{schmstad2006},
Theorem 6, under the additional assumption that the tail copula\vadjust{\goodbreak}
has continuous partial derivatives. As pointed out earlier, there is no
tail copula $\ltc\not\equiv0$ with this property.
We point out that in the case where $d=2$, it can be shown that
all tail copulas with continuous partial derivatives on the interior as
required in \eqref{fotc}
also have continuous partial derivatives on the axes, expect for the
origin. Therefore, in the two-dimensional case, our condition just
makes verification issues easier. Nevertheless,
in higher dimensions, the condition (as considered, e.g.,
Einmahl \textit{et~al.} \cite{einkraseg2011}) becomes more meaningful.
Note that a careful inspection of our proof in the \hyperref[appa]{Appendix} reveals
that a generalization to higher dimensions is not a trivial extension
(especially the proof of Lemma~\ref{lem:hadamard}).
Nevertheless, we are convinced that such an extension is possible.
Recently, B{\"{u}}cher and Volgushev \cite{buevol2011} derived a similar extension to the
$d$-dimensional case for the
usual empirical copula process. This program should be transferred to
the empirical tail
copula process and is deferred to future research.

\section{Multiplier bootstrap approximation} \label{sec:tcmult}
\subsection{Asymptotic theory} \label{subsec:asytail}
In this section we construct multiplier bootstrap approximations of the
Gaussian limit distributions $\Gb_{\ltctilde}$ and $\Gb_{\ltchat}$
specified in
\eqref{t4ss} and \eqref{t5ss}, respectively. Toward this end, let
$\xi_i$ be
independent identically distributed
positive random variables, independent of the $\vect{X}_i$, with mean
$\mu$ in
$(0,\infty)$ and finite variance $\tau^2$.
We first deal
with the case of known marginals and define a multiplier bootstrap analog
of \eqref{est:ltctilde} by
%
\begin{equation} \label{lamtil}
\ltctilde^\xi(\vect{x})=\frac{1}{k}\sum_{i=1}^n \frac{\xi
_i}{\bar\xi_n}\ind\biggl\{F_1(X_{i1})\leq\frac{kx_1}{n}, F_2(X_{i2})\leq
\frac{kx_2}{n}\biggr\},
\end{equation}
where $\bar{\xi}_n=n^{-1}\sum_{i=1}^n \xi_i$ denotes the mean of
$\xi_1, \ldots, \xi_n$. We have
%
\begin{equation} \label{eq:atildem}
\tilde{\alpha}^{m}_n(\vect{x})= \frac{\mu}{\tau}\frac{1}{\sqrt
{n}}\sum_{i=1}^n
\biggl( \frac{\xi_i}{\bar{\xi}_n}-1 \biggr)f_{n,\vect{x}}(U_i)=\frac{\mu
}{\tau}\sqrt{k} ( \ltctilde^\xi-\ltctilde),
\end{equation}
where the function $f_{n,\vect{x}}(U_i)$ is defined by
%
\begin{equation}\label{f_nx}
f_{n,\vect{x}} (\vect{U}_i) = \sqrt{\frac{n}{k}} \ind \{ U_{i1}
\leq kx_1/n, U_{i2} \leq kx_2/n \},
\end{equation}
and
\[
\vect{U}_i=(U_{i1},U_{i2}) ;\qquad U_{ip} = F_p(X_{ip}) \qquad\mbox{for } p=1,2.
\]
Throughout this paper, we use the notation
%
\begin{equation}\label{weakkos}
G_n \weakP G \qquad\mbox{in } \Db
\end{equation}
for \textit{conditional weak convergence in a metric space $(\Db,d)$}
in the sense of Kosorok \cite{kosorok2008}, page 19.
To be precise, \eqref{weakkos} holds\vadjust{\goodbreak} for some random variables
$G_n=G_n(\vect{X}_1,\ldots,\vect{X}_n,\xi_1,\ldots\xi_n),\break G\in
\Db$ if and only if
%
\begin{equation}\label{BL}
\sup_{h \in \mathit{BL}_1(\Db)} | \Exp_\xi h(G_n) - \Exp h(G)| \Pconv0
\end{equation}
and
%
\begin{equation}\label{am}
\Exp_\xi h(G_n)^*-\Exp_\xi h(G_n)_* \Pconv0 \qquad\mbox{for every } h
\in \mathit{BL}_1(\Db),
\end{equation}
where
\[
\mathit{BL}_1(\Db) = \{ f\dvt \Db\rightarrow\R \vert \Vert f\Vert_\infty\leq1, |f(\beta
)-f(\gamma)| \leq d(\beta,\gamma)\ \forall\gamma,\beta\in\Db \}
\]
denotes the set of all Lipschitz-continuous functions bounded by $1$.
The subscript $\xi$ in the expectations indicates conditional expectation
over the weights $\xi=(\xi_1,\ldots,\xi_n)$ given the data, and $h(G_n)^*$
and $ h(G_n)_*$ denote measurable majorants and minorants with respect
to the joint data, including the weights $\xi$. The condition (\ref
{BL}) is
motivated by the metrization of weak convergence by the bounded
Lipschitz metric (see, e.g., Theorem 1.12.4 in {V}an~der Vaart
\cite{vandervaart1998}). The
following result shows that the process \eqref{eq:atildem} provides a valid
bootstrap approximation of the process defined in \eqref{eq:atilde}.
\begin{theorem} \label{theo:atildem}
If $\ltc\not\equiv0$ and the second-order condition
\eqref{secornew} holds with $\sqrt{k}A(n/k)\rightarrow0 $,
$k=k(n)\rightarrow\infty$ and $k=\mathrm{o}(n)$
we have, as $n$ tends to infinity,
\[ \tilde{\alpha}_n^{m}=\frac{\mu}{\tau}\sqrt{k} (
\ltctilde^\xi-\ltctilde) \weakP\Gb_{\ltctilde} 
\]
in the metric space $\mathcal{B}_\infty(\bar{\R}^2_+)$.
\end{theorem}

Because Theorem~\ref{theo:atildem} states that we have weak
convergence of
$\tilde{\alpha}_n^{m}$ to $\Gb_{\ltctilde}$ conditional on the data
$U_i$, it
provides a bootstrap approximation of the empirical tail copula in the case
where the marginal distributions are known. To be precise, consider
$B\in\Nat$ independent replications of the random variables $\xi
_1,\ldots,\xi_n$
and denote them by $\xi_{1,b},\ldots,\xi_{n,b}$. Compute the statistics
$\tilde{\alpha}_{n,b}^{m}=\tilde{\alpha}_n^{m}(\xi_{1,b},\ldots
,\xi_{n,b})$ $(b=1,\ldots, B)$
and use the empirical distribution of $\tilde{\alpha
}_{n,1}^{m},\ldots,\tilde{\alpha}_{n,B}^{m}$
as an approximation for the limiting distribution of $G_{\ltctilde}$.

Because in most cases of practical interest there will be no information
about the marginals, Theorem~\ref{theo:atildem} cannot be used in many
statistical applications. We have developed two consistent bootstrap
approximation for the limiting distribution of the process \eqref{t5ss}
that do not require knowledge of the marginals. Intuitively, it is natural
to replace the unknown marginal distributions in \eqref{lamtil} by their
empirical counterparts, that is,
%
\begin{equation}
\label{lamdach}
\ltchat^{\xi,\cdot}(\vect{x})=\frac{1}{k}\sum_{i=1}^n \frac{\xi
_i}{\bar{\xi}_n} \ind\{ X_{i1}\leq F_{n1}^{-}(kx_1/n), X_{i2} \leq
F_{n2}^{-}(kx_2/n)\},
\end{equation}
which yields the process
\begin{eqnarray*}
\beta_n(\vect{x}) &=& \frac{\mu}{\tau}\sqrt{k} (\ltchat^{\xi
,\cdot}-\ltchat)
\\
&=& \frac{\mu}{\tau}\frac{1}{\sqrt{k}}\sum_{i=1}^n \biggl( \frac{\xi
_i}{\bar{\xi}_n}-1 \biggr)\ind\{ X_{i1}\leq F_{n1}^{-}(kx_1/n), X_{i2}
\leq F_{n2}^{-}(kx_2/n)\}.
\end{eqnarray*}
Unfortunately, this intuitive approach does not yield an approximation
for the distribution of the process $\Gb_{\ltchat}$, but only of $\Gb
_{\ltctilde}$.
\begin{theorem}\label{theo:betan}
Suppose that the assumptions of Theorem~\ref{theo:ltchat} hold.
Then we have, as $n$ tends to infinity,
\[
\beta_n=\frac{\mu}{\tau}\sqrt{k} ( \ltchat^{\xi,\cdot}-\ltchat
) \weakP\Gb_{\ltctilde} 
\]
in the metric space $\mathcal{B}_\infty(\bar{\R}^2_+)$.
\end{theorem}

Although Theorem~\ref{theo:betan} provides a negative result and shows
that the distribution of $\beta_n$ cannot be used for approximating
the limiting law $\Gb_{\ltchat}$, it turns out to be essential for
our first consistent
multiplier bootstrap method. To be precise, we note that the
distribution of
$\beta_n$ can be calculated from the data without knowing the
marginal distributions. Consequently, we obtain an approximation for
the unknown distribution of the process $\Gb_{\ltctilde}$. To get an
approximation of $\Gb_{\ltchat}$, we follow R\'{e}millard and Scaillet \cite{remiscai2009} and estimate
the derivatives of the tail copula by
\[
\widehat{\partial_p \ltc}(\vect{x}) :=
\cases{
\displaystyle\frac{\ltchat(\vect{x}+h\vect{e}_p)-\ltchat(\vect{x}-h\vect
{e}_p)}{2h}, &\quad $\infty>x_p\geq h$,\vspace*{2pt}
\cr
\widehat{\partial_p \ltc}(\vect{x} + (h-x_p) \vect{e}_p), &\quad $x_p <
h$,\vspace*{2pt}
\cr
0, &\quad $x_p=\infty$,}%
\]
where $\vect{e}_p$ denotes the $p$th unit vector ($p=1,2$) and
$h \sim k^{-1/2}$ tends to $0$ with increasing sample size. In the
\hyperref[appa]{Appendix} (see the proof of the following theorem in
Appendix~\ref{appa}), we
show that these
estimates
are consistent, and thus we define the process
%
\begin{equation} \label{alphapdm}
\alpha^{\mathit{pdm}}_n(\vect{x}) = \beta_n(\vect{x}) - \widehat{\partial
_1 \ltc}(\vect{x})\beta_{n}(x_1,\infty) - \widehat{\partial_2
\ltc}(\vect{x})\beta_{n}(\infty,x_2).
\end{equation}
Note that $\alpha_n^{\mathit{pdm}}$ depends only on the data and the multipliers
$\xi_1,\ldots,\xi_n$. Consequently, a bootstrap sample can be
readily generated
as described in the previous paragraph; in what follows, we call this
method the \textit{partial derivatives
multiplier bootstrap ($\mathit{pdm}$-bootstrap)}. Our next result shows
that the $\mathit{pdm}$ bootstrap provides a valid approximation for the
distribution of the process
$\Gb_{\ltchat}$.

\begin{theorem} \label{theo:alphanpdm}
Under the assumptions of Theorem~\ref{theo:ltchat}, we have
\[ \alpha^{\mathit{pdm}}_n \weakP\Gb_{\ltchat}
\]
in the metric space $\mathcal{B}_\infty(\bar{\R}^2_+) $.
\end{theorem}

It turns out that there is an alternative valid multiplier bootstrap
procedure in
the case of unknown marginal distributions, which is attractive because it
avoids the problem of estimating the partial derivatives of the lower
tail copula.
This method introduces multiplier random variables not only in the
two-dimensional
distribution function, but also in the inner estimators of the
marginals. To be precise,
define
\begin{eqnarray*}
F_n^\xi(\vect{x}) &=& \frac{1}{n}\sum_{i=1}^n \frac{\xi_i}{\bar
\xi_n} \ind\{ X_{i1}\leq x_1, X_{i2}\leq x_2\},
\\
F_{nj}^\xi(x_j) &=& \frac{1}{n}\sum_{i=1}^n \frac{\xi_i}{\bar\xi
_n} \ind\{ X_{ij}\leq x_j\},\qquad j=1,2,
\\
C_n^{\xi,\xi}(\vect{u}) &=& F_n^\xi(F_{n1}^{\xi-}(u_1), F_{n2}^{\xi-}(u_2))
\end{eqnarray*}
and consider the process
%
\begin{equation} \label{ltchat2xi}
\ltchat^{\xi,\xi}(\vect{x})  := \frac{n}{k} C_n^{\xi,\xi}
\biggl(\frac{k}{n}\vect{x} \biggr)
=\frac{1}{k}\sum_{i=1}^n \frac{\xi_i}{\bar{\xi}_n} \ind\{
X_{i1}\leq{F}^{\xi-}_{n1}(kx_1/n), X_{i2} \leq{F}^{\xi
-}_{n2}(kx_2/n) \}.
\end{equation}
We call this bootstrap method the \textit{direct multiplier
bootstrap ($\mathit{dm}$ bootstrap)}.

\begin{theorem}\label{theo:alphandm}
Under the assumptions of Theorem~\ref{theo:ltchat}, we have
%
\begin{equation}\label{eq:alphandm}
\alpha_n^{\mathit{dm}} = \frac{\mu}{\tau} \sqrt{k} (\ltchat
^{\xi,\xi}-\ltchat )
\weakP\Gb_{\ltchat}\qquad \mbox{in } \mathcal{B}_\infty(\bar{\R
}^2_+).
\end{equation}
\end{theorem}

\begin{rem}As pointed out by a referee, an alternative multiplier
bootstrap could be obtained by multiplying each summand with $\xi
_i-1$, where $\xi_1,\ldots,\xi_n$ are i.i.d. with
$E[\xi_1]=\operatorname{Var}(\xi_i)=1$ (see Kosorok \cite{kosorok2008}, Section 11.4.2).
\end{rem}

\subsection{Finite-sample results} \label{subsec:fstail}
In this section we present a small comparison of the finite-sample
properties of the two bootstrap approximations given in this section.
We also study the impact of the choice of the parameter $k$ on the
properties of the estimates and the bootstrap
procedure.
For the sake of brevity, we only consider data generated form the
Clayton copula with a coefficient of lower tail dependence $\lambda_L=0.25$.
The Clayton copula, defined by
%
\begin{equation}\label{claytoncop}
C(\vect{u};\theta)= (u_1^{-\theta}+u_2^{-\theta}-1 )^{-1/\theta},\qquad
\theta>0,
\end{equation}
is widely used for modeling of negative tail-dependent data.
Its lower tail copula is given by
\[
\ltc(\vect{x})= (x_1^{-\theta}+x_2^{-\theta} )^{-1/\theta}.
\]
Tables~\ref{tab:covtc1} and~\ref{tab:covtc2} show the accuracy of the
bootstrap approximation of the covariances of the
limiting  variable $\Gb_{\ltchat}$.

\begin{figure}

\includegraphics{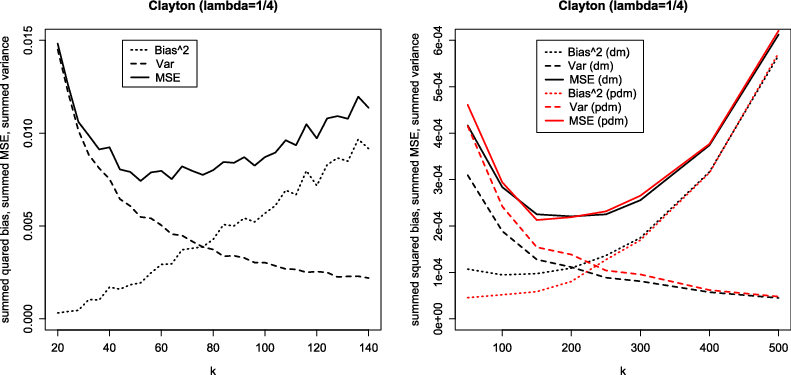}

\caption{\emph{Left: Averaged MSE, variance, and squared bias for the
estimation of
$\ltc(\mathrm{e}^{\mathrm{i}\ell\uppi/4})$ ($\ell=1,2,3$) by its empirical counterpart
$\ltchat(\mathrm{e}^{\mathrm{i}\uppi/4})$ against the parameter $k$. Right:
Averaged MSE, variance, and squared bias for the bootstrap estimation
of the covariances of $\Gb_{\ltchat}$.}}\label{pic:cov}
\end{figure}

We chose three points on the unit circle, $\{ \mathrm{e}^{\mathrm{i}\varphi}, \varphi
={\ell\uppi/8} \mbox{ with } \ell=1,2,3\}$,
and present in the first four columns of Table~\ref{tab:covtc1}, the
true covariances of the limiting process $\Gb_{\ltchat}$. The
remaining columns show the simulated
covariances of the process $\alpha_n$ on the basis of $5 \cdot10^5$
simulation runs, with a sample size of $n=1000$ and the parameter $k$
chosen as $50$.
This choice is motivated by the left panel of Figure~\ref{pic:cov},
which plots the sum of the
squared bias, the variance, and the mean squared error (MSE) of the
estimators $\ltchat(\mathrm{e}^{\mathrm{i}\ell\uppi/4})$ for $\ltc(\mathrm{e}^{\mathrm{i}\ell \uppi/4})$
($\ell=1,2,3$). The MSE is minimized for values of $k$
in a neighborhood of the point $50$. Note also
that the literature provides several data-adaptive proposals for the
choice of the parameter $k$ (see, e.g., Drees and Kaufmann \cite{dreekauf1998} oder
Gomes and Oliveira~\cite{gomeoliv2001}) in the univariate context.

\begin{table}
\tablewidth=\textwidth
\tabcolsep=0pt
  \caption{\textit{Left: True covariances of $\Gb_{\ltchat}$ for the Clayton
copula with $\lambda_{L}=0.25$. Right: Sample covariances of the
empirical tail copula process $\alpha_n$ with sample size
$n=1000$ and parameter $k=50$}}\label{tab:covtc1}
\begin{tabular*}{\textwidth}{@{\extracolsep{\fill}}lllllll@{}}
\hline
& \multicolumn{3}{l}{True} & \multicolumn{3}{l}{$\alpha_n$}
\\[-5pt]
&\multicolumn{3}{l}{\hrulefill} & \multicolumn{3}{l@{}}{\hrulefill}
\\
& $\frac{\uppi}{8}$ & $2\frac{\uppi}{8}$ & $ 3\frac{\uppi}{8}$ & $\frac
{\uppi}{8}$ & $2\frac{\uppi}{8}$ & $ 3\frac{\uppi}{8}$
\\
\hline
$\frac{\uppi}{8}$ & 0.0874 & 0.0754 & 0.0516 & 0.0889 & 0.0737 & 0.0476
\\[2pt]
2$\frac{\uppi}{8}$ & & 0.1160 & 0.0754 & & 0.1218 & 0.0741
\\[2pt]
3$\frac{\uppi}{8}$ & & & 0.0874 & & & 0.0892
\\
\hline
\end{tabular*}
\end{table}

\begin{table}[b]
\tablewidth=\textwidth
\tabcolsep=0pt
\caption{\textit{Averaged sample covariances (rows 3--5) and $\mathrm{MSE} \times10^4$
(rows 6--8) of the bootstrap approximations $\alpha_n^{\mathit{pdm}}$, $\alpha
_n^{\mathit{dm}}$ and $\alpha_n^{\mathit{res}}$ of $\Gb_{\ltchat}$ under the
conditions of Table \protect\hyperref[tab:covtc1]{\emph{1}}}} \label{tab:covtc2}
\begin{tabular*}{\textwidth}{@{\extracolsep{\fill}}llllllllll@{}} \hline
&\multicolumn{3}{l}{$\alpha_n^{\mathit{pdm}}$} &\multicolumn{3}{l}{$\alpha_n^{\mathit{dm}}$}
& \multicolumn{3}{l}{$\alpha_n^{\mathit{res}}$}
\\[-5pt]
&\multicolumn{3}{l}{\hrulefill}&\multicolumn{3}{l}{\hrulefill}
 &\multicolumn{3}{l}{\hrulefill}
 \\
& $\frac{\uppi}{8}$ & $2\frac{\uppi}{8}$ & $ 3\frac{\uppi}{8}$ & $\frac
{\uppi}{8}$ & $2\frac{\uppi}{8}$ & $ 3\frac{\uppi}{8}$ & $\frac{\uppi
}{8}$ & $2\frac{\uppi}{8}$ & $ 3\frac{\uppi}{8}$
\\
\hline
$\frac{\uppi}{8}$ & 0.094 & 0.072 & 0.046 & 0.100 & 0.071 & 0.045 &
0.100 & 0.070 & 0.043
\\[2pt]
$2\frac{\uppi}{8}$ & & 0.130 & 0.072 & & 0.136 & 0.707 & & 0.136 &
0.070
\\[2pt]
$3\frac{\uppi}{8}$ & & & 0.094 & & & 0.099 & & & 0.099
\\[5pt]
$\frac{\uppi}{8}$ & 3.67 & 4.68 & 3.65 & 3.86 & 3.49 & 2.72 & 4.21 &
3.85 & 3.21
\\[2pt]
$2\frac{\uppi}{8}$ & & 8.11 & 4.87 & & 8.89 & 3.25 & & 8.73 & 3.64
\\ [2pt]
$3\frac{\uppi}{8}$ & & & 3.70 & & & 3.77 & & & 3.90
\\
\hline
\end{tabular*}
\end{table}

The data in Table~\ref{tab:covtc1} serve as a benchmark for the
multiplier bootstrap approximations of the covariances stated in Table
\ref{tab:covtc2},
where we investigate the quality of the approximation by various
bootstrap methods.
The distribution of the multipliers in the $\mathit{dm}$ and $\mathit{pdm}$ bootstrap procedures
was chosen according to B{\"{u}}cher and Dette \cite{buecdett2010} as $\Prob(\xi=0)=\break\Prob
(\xi=2)=0.5$, such that $\mu=\allowbreak\tau=1$.
For the sake of completeness, we also investigate
the resampling bootstrap considered in Peng and Qi \cite{pengqi2008}, hereinafter
denoted by $\alpha_n^{\mathit{res}}$. The estimated covariances given in
the first part (rows 3--5) of Table~\ref{tab:covtc2} were calculated
by $1000$ simulation
runs, where in each run the covariance is estimated on the basis of
$B=500$ bootstrap replications.
The second part (rows 6--8) of Table~\ref{tab:covtc2} shows the
corresponding MSE.

As the figures show, all bootstrap procedures yield approximations of
comparable quality.
Considering only the bias in Table~\ref{tab:covtc2}, the $\mathit{pdm}$
bootstrap demonstrates slight advantages in all cases,
whereas there are basically no differences between the $\mathit{dm}$ and
the resampling bootstraps.
A comparison of the MSE in Table~\ref{tab:covtc2} shows that the
$\mathit{pdm}$ bootstrap has the best performance
on the diagonal. On the other hand, it yields a less accurate
approximation for the off-diagonal covariances, for which the \textit
{$\mathit{dm}$} bootstrap yields the best results.

The right panel of Figure~\ref{pic:cov} illustrates the sensitivity of
the accuracy of the estimators for the covariances with respect to the
choice of $k$. For this purpose, we calculated
the sum of the MSE values given in Table~\ref{tab:covtc2} (as well as
the variance and squared bias) for various choices of~$k$.
As can be seen, the best choices for $k$ lie in an interval of
approximate length $100$ around the center $k=200$. Compared with the
``best'' value $k=50$, for estimating $\ltc$, the optimal values for
estimating the covariances of $\Gb_{\ltchat}$ are approximately four
times larger for both the $\mathit{pdm}$ bootstrap and the $\mathit{dm}$ bootstrap.
This increase may be explained by the fact that the large bias of
$\ltchat(\vect{x}), \ltchat^{\xi,\cdot}(\vect{x})$ and $\ltchat
^{\xi,\xi}(\vect{x})$ for estimating $\ltc(\vect{x})$ cancels out
if the difference $\ltchat^{\xi,\cdot}(\vect{x})- \ltchat(\vect
{x})$ or $\ltchat^{\xi,\xi}(\vect{x})- \ltchat(\vect{x})$ is
calculated. As a result, we may choose larger values of $k$, resulting
in a notable decay of the variance.

These findings have ambiguous consequences. If we are interested only
in the covariances of the limiting variable $\Gb_{\ltchat}$, then a
larger value of $k$ for the bootstrap is\vspace*{1pt} advisable. However, because
the bootstrap is not able to capture the true bias of the empirical
tail copula process, some care is needed if we are interested in
approximations of the whole distribution of $\sqrt{k}(\ltchat-\ltc)$
(as is required in Section~\ref{sec:tcapplications}). In this case a
careful choice of $k$ for estimating the tail copula becomes even more
important; we are confronted with the strong requirement of a small
bias for this estimator.
Finally, a comparison of the variance and the bias of the two bootstrap
procedures investigated in Figure~\ref{pic:cov} reveals that the $\mathit{pdm}$
bootstrap has a smaller bias, but a slightly larger variance, than the
$\mathit{dm}$ bootstrap. On the other hand, the differences with respect to the
MSE are nearly undetectable.
%
%
%
%
%
\section{Statistical applications} \label{sec:tcapplications}
In this section we investigate several statistical applications of the
multiplier bootstrap. In particular, we discuss the problem of comparing
lower tail copulas from different samples, the problem of constructing
confidence intervals, and the problem of testing for a parametric form
of the lower tail copula.

\subsection{Testing for equality between two tail copulas}\label
{subsec:equality}

Let $\vect{X}_1,\ldots,\vect{X}_{n_1}$ and $\vect{Y}_1,\ldots
,\vect{Y}_{n_2}$ denote
two independent samples of i.i.d.
random variables (we relax the assumption of
independence between the samples later) with continuous cumulative
distribution function
$F=C(F_1,F_2)$ and $H=D(H_1,H_2)$, respectively. We assume that
for both distributions, the corresponding lower tail copulas, say
$\ltci{X}$ and $\ltci{Y}$,
exist and do not vanish and the tails of the corresponding copulas
converge to $\ltc$ at the rate specified in \eqref{secornew}. We are
interested in a test for the hypothesis
%
\begin{equation} \label{h0comp}
\mathcal{H}_0\dvt \ltci{X}\equiv\ltci{Y} \quad\mbox{vs.}\quad \mathcal{H}_1\dvt
\ltci{X}\not\equiv\ltci{Y}.
\end{equation}
Given the homogeneity of tail copulas, we have
$\ltc(t\vect{x})=t\ltc(\vect{x})$ for all $t>0, \vect{x}\in
[0,\infty)^2$,
and the hypotheses are equivalent to
\[
\mathcal{H}_0\dvt \varrho(\ltci{X},\ltci{Y})=0 \quad\mbox{vs.}\quad \mathcal
{H}_1\dvt \varrho(\ltci{X},\ltci{Y})>0,
\]
where the distance $\varrho$ is defined by
%
\begin{eqnarray} \label{dist}
\varrho(\ltci{X},\ltci{Y})&:=& \int_0^{\uppi/2} \bigl( \ltci{X}(\cos
\varphi,\sin\varphi)- \ltci{Y}(\cos\varphi,\sin\varphi) \bigr)^2
\,\mathrm{d}\varphi\nonumber
\\[-8pt]
\\[-8pt]
&\hspace{2.7pt}=& \int_0^{\uppi/2} \bigl( \ltciein{X} (\varphi)- \ltciein{Y}(\varphi)
\bigr)^2 \,\mathrm{d}\varphi\nonumber
\end{eqnarray}
and we use the notation $\ltciein{X} (\varphi)= \ltci{X}(\cos
\varphi,\sin\varphi), \ltciein{Y} (\varphi)= \ltci{Y}(\cos
\varphi,\sin\varphi)$.
{Note that the integration in \eqref{dist} over the unit circle with
respect to the Euclidian norm is rather a matter of taste. An
integration along the unit circle with respect to the sup-norm on $\R
^2$ (i.e., along $\{\vect{x}\in[0,\infty)^2\dvt \max(x_1,x_2)=1\}$)
would be possible as well.}

We propose basing the test
for the hypothesis \eqref{h0comp}
on the distance between the empirical tail copulas
and define
\[
\mathcal{S}_n=\frac{k_1k_2}{k_1+k_2} \varrho(\ltcihat{X},\ltcihat{Y})
=\frac{k_1k_2}{k_1+k_2}\int_0^{\uppi/2}\bigl(\ltcihatein{X} (\varphi)
-\ltcihatein{Y} (\varphi)\bigr)^2 \,\mathrm{d}\varphi,
\]
where $\ltcihatein{X} (\varphi) = \ltcihat{X} (\cos(\varphi), \sin
(\varphi))$, $\ltcihatein{Y}=
\ltcihat{Y} (\cos(\varphi), \sin(\varphi))$ denote the empirical
tail copulas $\ltcihat{X}$
and $\ltcihat{Y}$ with corresponding parameters $k_1$ and $k_2$, satisfying
\[
k_p\rightarrow\infty, \qquad k_p=\mathrm{o}(n_p) \qquad(p=1,2) \quad\mbox{and}\quad
k_1/(k_1+k_2)\rightarrow\lambda\in(0,1) .
\]
We assume that the tail copulas $\ltci{X}$ and $\ltci{Y}$ satisfy a
second-order condition as in \eqref{secornew} (with $A$ replaced by
$A_p$), and that $k_p$ is chosen appropriately, that is, $\sqrt{k_p}
A_p(k_p/n_p)=\mathrm{o}(1)$.
Under the null hypothesis \eqref{h0comp} of equality between the tail copulas,
we have $\mathcal{S}_n=\mathcal{T}_n$ with
\[
\mathcal{T}_n = \int_0^{\uppi/2} \Ec_n^2 (\cos\varphi, \sin\varphi
) \,\mathrm{d} \varphi,
\]
where
\[
\Ec_n( \vect{x} ) = \sqrt{\frac{k_2}{k_1+k_2}} \sqrt{k_1}\bigl(\ltcihat
{X}(\vect{x})-\ltci{X}(\vect{x} )\bigr) -
\sqrt{\frac{k_1}{k_1+k_2}} \sqrt{k_2} \bigl(\ltcihat{Y}( \vect{x}
)-\ltci{Y}( \vect{x})\bigr).
\]
Because the two samples $X$ and $Y$ are independent, we obtain, independent
of the hypotheses, that
%
\begin{equation}\label{Ec}
\Ec_{n} \weak\sqrt{1-\lambda} \Gb_{\ltcihat{X}} - \sqrt{\lambda
} \Gb_{\ltcihat{Y}} =: \Ec_{\cdot}
\end{equation}
in the metric space $\mathcal{B}_\infty(\bar{\R}^2_+)$, where the
stochastically independent two-dimensional
centered Gaussian fields $\Gb_{\ltcihat{X}}$ and $\Gb_{\ltcihat
{Y}}$ are as defined in (\ref{gf}). This yields, by the continuous
mapping theorem,
\[
\mathcal{T}_n\weak\int_{0}^{\uppi/2} \Ec^2(\cos\varphi,\sin
\varphi)\,\mathrm{d}\varphi=:\mathcal{T}
\]
under both the null hypothesis and the alternative. Note that
$\varrho(\ltcihat{X},\ltcihat{Y})\stackrel{\Prob}{\rightarrow
}\varrho(\ltci{X},\ltci{Y})$,
which vanishes if and only if the null hypothesis \eqref{h0comp} is
satisfied. Therefore,
we can conclude that
%
\begin{equation}\label{S_n}
\mathcal{S}_n \weak_{\mathcal{H}_0} \mathcal{T}, \qquad \mathcal{S}_n
\stackrel{\Prob}{\rightarrow}_{\mathcal{H}_1} \infty,
\end{equation}
which shows that a test that rejects the null hypothesis \eqref
{h0comp} for large
values of $\mathcal{T}_n$ is consistent. {Note that the latter
convergence depends crucially on the assumption that
$\sqrt{k}A(n/k)\rightarrow0$. If this assumption does not hold, then
a large value of $\mathcal{S}_n$ could
occur even under the null hypothesis, because of the biasedness of
$\ltchat$. Thus, as discussed at the end of Section~\ref
{subsec:fstail}, the choice of a small $k$ corresponding to a small
bias is of prime importance.}

To determine critical values for the test, we approximate the limiting
distribution $\mathcal{T}$ by the multiplier bootstrap proposed in Section
\ref{sec:tcmult}. For this purpose, we consider the $\mathit{pdm}$ bootstrap
(the extension to the $\mathit{dm}$ bootstrap is straightforward) using the
definition in
equation (\ref{eq:alphandm}) and denote, for any $b\in\{1,\ldots,B\}
$, $\xi_{1,b},\ldots,\xi_{n_1,b}, \zeta_{1,b},\ldots,\zeta
_{n_2,b}$ i.i.d. nonnegative random variables with mean $\mu_1$
(resp., $\mu_2$) and variance $\tau_1^2$ (resp., $\tau_2^2$). For each
$b$ and both samples, we compute the bootstrap statistics as given in
\eqref{alphapdm}, that is,
\begin{eqnarray*}
\alpha^{\mathit{pdm}}_{X,n_1,b}(\vect{x}) &=& \beta_{X,n_1,b}(\vect{x}) -
\widehat{\partial_1 \ltci{X}}(\vect{x})\beta_{X,n_1,b}(x_1,\infty
) - \widehat{\partial_2 \ltci{X}}(\vect{x})\beta_{X,n_1,b}(\infty
,x_2),
\\
\alpha^{\mathit{pdm}}_{Y,n_2,b}(\vect{x}) & =& \beta_{Y,n_2,b}(\vect{x}) -
\widehat{\partial_1 \ltci{Y}}(\vect{x})\beta_{Y,n_2,b}(x_1,\infty
) - \widehat{\partial_2 \ltci{Y}}(\vect{x})\beta_{Y,n_2,b}(\infty,x_2),
\end{eqnarray*}
where
\begin{eqnarray*}
\beta_{X,n_1,b}(\vect{x}) &=& \frac{\mu_1}{\tau_1}\frac{1}{\sqrt
{k_1}}\sum_{i=1}^{n_1} \biggl( \frac{\xi_{i,b}}{\bar{\xi}_{\cdot
,b_{n_1}}}-1 \biggr)\ind\{ X_{i1}\leq F_{n_11}^-(k_1x_1/n_1), X_{i2}\leq
F_{n_12}^-(k_1x_2/n_1)\},
\\
\beta_{Y,n_2,b}(\vect{x}) &=& \frac{\mu_2}{\tau_2}\frac{1}{\sqrt
{k_2}}\sum_{i=1}^{n_2} \biggl( \frac{\zeta_{i,b}}{\bar{\zeta}_{\cdot
,b_{n_2}}}-1 \biggr)\ind\{ Y_{i1}\leq H_{n_21}^-(k_2x_1/n_2), Y_{i2}\leq
H_{n_22}^-(k_2x_2/n_2)\},
\end{eqnarray*}
and $\widehat{\partial_{p} \ltci{X}}$ and $\widehat{\partial_{p}
\ltci{Y}}$
are the corresponding\vspace*{1pt} estimates of the partial
derivatives\break
($p=1,2$). For all $\vect{x}\in\bar\R_+^2$ and all $b\in\{1,\ldots
,B\}$, define
\begin{eqnarray*}
\hat\Ec_{n}^{(\mathit{pdm}, b)} (\vect{x}) &:=& \sqrt{\frac{k_2}{k_1+k_2}}
\alpha^{\mathit{pdm}}_{X,n_1,b}(\vect{x}) - \sqrt{\frac{k_1}{k_1+k_2}}
\alpha^{\mathit{pdm}}_{Y,n_2,b}(\vect{x}),
\\
\hat{\mathcal{T}}_{n}^{(\mathit{pdm},b)} & :=& \int_{0}^{\uppi/2} \bigl\{ \hat\Ec
_{n}^{(\mathit{pdm},b)} ( \cos\varphi, \sin\varphi) \bigr\}^2 \,\mathrm{d} \varphi.
\end{eqnarray*}
By Theorem~\ref{theo:alphanpdm} and Theorem 10.8 in Kosorok
\cite{kosorok2008}, it follows that for every $b\in\{1,\ldots,B\}$,
\[
\hat{\mathcal{T}}_{n}^{(\mathit{pdm},b)} \weakP\mathcal{T}^{(b)},
\]
where $\mathcal{T}^{(b)}$ is an independent copy of $\mathcal{T}$.
(Note that
we consider the processes $\hat\Ec_{n}^{(\mathit{pdm},b)}$ in the Banach space
$\ell^\infty([0,1]^2)$.)
Thus, from (\ref{S_n}), we obtain a consistent asymptotic level
$\alpha$
test for the null hypothesis \eqref{h0comp} by rejecting $\mathcal
{H}_0$ for
large values of $\mathcal{S}_{n}$, that is,
%
\begin{equation} \label{testcomp}
\mathcal{S}_{n} > q_{1-\alpha}^{\mathit{pdm}},
\end{equation}
where $q_{1-\alpha}^{\mathit{pdm}}$ denotes the $(1-\alpha)$-quantile of the
c.d.f. \mbox{$
K_{n}^{\mathit{pdm}}(s) = B^{-1}\sum_{b=1}^B \ind\{ \hat{\mathcal
{T}}_{n}^{(\mathit{pdm},b)} \leq s \}.
$}

So far, we have focused our discussion on the case of two independent
samples. It is easy to check that our\vadjust{\goodbreak} methodology also applies in cases
of paired observations, that is, $\vect{X}_i$ is not independent of
$\vect{Y}_i$,
but $n_1=n_2=n$. In that case we must set $\zeta_{i,b}=\xi_{i,b}$ for
all $i=1,\ldots,n$ and $b=1,\ldots,B$. To see this, set
$\vect{Z}_i=(\vect{X}_{i1},\vect{X}_{i2},\vect{Y}_{i1},\vect
{Y}_{i2})$ and
denote the (empirical) copula of $\vect{Z}_i$ by $(\mathcal{C}_n$)
$\mathcal{C}$.
Clearly,
\begin{eqnarray*}
C(u_1,u_2)&=&\mathcal{C}(u_1,u_2,1,1), \qquad D(v_1,v_2)=\mathcal
{C}(1,1,v_1,v_2),
\\
C_n(u_1,u_2)&=&\mathcal{C}_n(u_1,u_2,1,1), \qquad D_n(v_1,v_2)=\mathcal
{C}_n(1,1,v_1,v_2).
\end{eqnarray*}
If we set $
\ltci{Z}(\vect{x},\vect{y})=\lim_{t\rightarrow\infty}t \mathcal
{C}(\vect{x}/t,\vect{y}/t)$, $ \ltcihat{Z}(\vect{x},\vect
{y})=\frac{n}{k}\mathcal{C}_n(\frac{n\vect{x}}{k},\frac{n\vect{y}}{k}),
$
then we obtain
\begin{eqnarray*}
\ltci{X}(\vect{x})&=&\ltci{Z}(\vect{x},\infty,\infty), \qquad
\ltci{Y}(\vect{y})=\ltci{Z}(\infty,\infty,\vect{y}),
\\
\ltcihat{X}(\vect{x})&=&\ltcihat{Z}(\vect{x},\infty,\infty),\qquad
\ltcihat{Y}(\vect{y})=\ltcihat{Z}(\infty,\infty,\vect{y}).
\end{eqnarray*}
Under a second-order condition on the joint tail copula $\ltci{Z}$,
the asymptotic properties
of the process $\ltcihat{Z}$ can be derived along similar lines as
before; we omit the details for the sake of brevity. As the only
difference from the preceding discussion, note that the occurring
limiting fields $\Gb_{\ltcihat{X}}$ and $\Gb_{\ltcihat{Y}}$ are no
longer independent.
Because the asymptotic behavior of the multiplier bootstrap
approximations can be shown to reflect this dependence, we still obtain
consistency of the test; we again omit the details.

\begin{table}[b]
\tablewidth=\textwidth
\tabcolsep=0pt
\caption{\textit{Simulated rejection probabilities of the
bootstrap tests defined in~\protect\hyperref[testcomp]{\emph{(4.5)}} for
the hypothesis~\protect\hyperref[h0comp]{\emph{(4.1)}}}}\label{tab:pdmequal}
\begin{tabular*}{\textwidth}{@{\extracolsep{\fill}}lllllllll@{}}
\hline
& & & \multicolumn{3}{l}{$\mathit{pdm}$} & \multicolumn{3}{l}{$\mathit{dm}$}
\\[-5pt]
& & &\multicolumn{3}{l}{\hrulefill} & \multicolumn{3}{l@{}}{\hrulefill}
\\
$k$ & $\lambda_{L,X}$ & $\lambda_{L,Y}$ & $\alpha=0.15$ & $\alpha
=0.1$ & $\alpha=0.05$ & $\alpha=0.15$ & $\alpha=0.1$ & $\alpha
=0.05$
\\
\hline
\phantom{0}50 &0.25 & 0.25 & 0.143 & 0.098 & 0.054 & 0.125 & 0.091 & 0.052
\\
& 0.5 & 0.5 & 0.140 & 0.099 & 0.047 & 0.108 & 0.069 & 0.036
\\
& 0.75 & 0.75 & 0.117 & 0.078 & 0.029 & 0.068 & 0.051 & 0.023
\\[3pt]
& 0.25 & 0.5 & 0.764 & 0.706 & 0.605 & 0.713 & 0.643 & 0.529
\\
& 0.5 & 0.75 & 0.896 & 0.856 & 0.783 & 0.869 & 0.822 & 0.713
\\
& 0.25 & 0.75 & 1 & 1 & 1 & 0.999 & 0.999 & 0.997
\\[6pt]
200 &0.25 & 0.25 & 0.145 & 0.107 & 0.052 & 0.125 & 0.084 & 0.044
\\
& 0.5 & 0.5 & 0.128 & 0.083 & 0.037 & 0.140 & 0.097 & 0.051
\\
& 0.75 & 0.75 & 0.141 & 0.092 & 0.041 & 0.103 & 0.068 & 0.035
\\[3pt]
& 0.25 & 0.5 & 0.991 & 0.978 & 0.948 & 0.979 & 0.971 & 0.950
\\
& 0.5 & 0.75 & 1 & 1 & 1 & 1 & 1 & 1
\\
& 0.25 & 0.75 & 1 & 1 & 1 & 1 & 1 & 1
\\
\hline
\end{tabular*}
\end{table}

To investigate the finite-sample property, we consider two
independent samples of i.i.d. random variables with
Clayton copula (see \eqref{claytoncop}) with a coefficient of lower
tail dependence
$\lambda_L$ varying in the set $\{0.25, 0.5, 0.75\}$.

Table~\ref{tab:pdmequal}
presents the simulated
rejection probabilities of the $\mathit{pdm}$ and $\mathit{dm}$ bootstrap tests defined in
\eqref{testcomp} for various nominal levels on the basis of $1000$ simulation
runs. The sample size was $n_1=n_2=n=1000$ and, $B=500$ bootstrap
replications with $\mathcal{U}(\{0,2\})$ multipliers (i.e., $\Prob
(\xi=0)=\Prob(\xi=2)=0.5$,
such that $\mu=\tau=1$) were used. The parameter $k$ was chosen as
either $k=50$ or $k=200$ as suggested by the discussion in the
preceding paragraph.

We observe that the nominal level is well approximated by the $\mathit{pdm}$ bootstrap
if the coefficient of tail dependence is not too large. For a larger
coefficient, the test tends to be
conservative.
Of note, the approximation of the nominal level is rather robust with
respect to the choice of~$k$.
A~comparison of the performance of the two bootstrap procedures shows
that the $\mathit{dm}$ bootstrap test is slightly more conservative,
and that this effect increases with the coefficient of tail dependence.

The alternative
of different lower tail copulas is detected with reasonable power. Both tests
yield rather similar results, with a slight advantage for the $\mathit{pdm}$ bootstrap.
Evaluation of the impact of the choice of the parameter $k$ under the
alternative shows
some advantages for $k=200$. This
again may be explained by the fact that bias terms cancel out if the
difference $\ltcihat{X}-\ltcihat{Y}$ is calculated.

%
%
%
%
\subsection{Bootstrap approximation of a minimum distance estimate
and a computationally efficient goodness-of-fit test}

In this section we estimate the tail copula of $\vect{X}$
under the additional assumption that it is an element of some parametric
class, say
$ \mathcal{L} = \{ \ltc(\cdot; \theta) \vert \theta\in\Theta\}$.
Estimation of parametric classes of tail copulas and stable tail
dependence functions was recently investigated by de~Haan \textit{et~al.} \cite{dehnevpen2008}
and Einmahl \textit{et~al.} \cite{einkraseg2008}, who proposed a censored likelihood
estimator and a
moment-based estimator, respectively. Here we investigate
a different, based on the minimum distance method. To be
precise, let $\ltc$ denote an arbitrary lower tail copula and $\ltc(
\cdot; \theta)$
denote an element in the parametric class $\mathcal{L} $, and
consider
the parameter corresponding to the best approximation by the distance
$\varrho$ defined
in \eqref{dist}
%
\begin{equation} \label{bestpar}
\theta_B = T(\ltc) = \arg\min_{\theta\in\Theta} \varrho(\ltc
,\ltc( \cdot; \theta)),
\end{equation}
where $\varrho$ is as defined in \eqref{dist}. We call $\thetamd
=T(\ltchat)$ a minimum distance estimator for $\theta$,
where $\ltchat$ is the empirical lower tail copula defined in (\ref{empl}).
Note that $\theta_B$ is the
``true'' parameter if and only if the null hypothesis is satisfied.

Throughout this section, we let $\vect{X}_1,\ldots,\vect{X}_n$ denote
i.i.d. bivariate
random variables with c.d.f. $F=C(F_1,F_2)$ and existing
lower tail copula $\ltc$. (For a proof of the following result, see
Section~4.5 in B{\"{u}}cher \cite{buecher2011}.)
%
%
%
\begin{theorem}\label{theo:thetamdasy}
Suppose that the standard conditions of minimum distance estimation are
satisfied. (For a precise formulation of these conditions, see
B{\"{u}}cher \cite{buecher2011}, pages 89ff.)
If the true tail copula
$\ltc$ satisfies the first-order condition \eqref{fotc} of Theorem
\ref{theo:ltchat}, and if the second-order condition \eqref{secornew}
holds with
$\sqrt{k}A(n/k)\rightarrow0 $, where $k=k(n)\rightarrow\infty$ and
$k=\mathrm{o}(n)$, then the minimum distance
estimator $\thetamd$ is consistent for the  parameter $\theta
_B$
corresponding to the best approximation with respect to the
distance $\varrho$.
Moreover,
\begin{eqnarray*}
\Theta^{\mathrm{MD}}_n &:=& \sqrt{k}( \thetamd-\theta_B) = \sqrt{k}\int
\gamma_{\theta_B}(\varphi) \bigl( \ltchatein(\varphi)-\ltcein(\varphi
) \bigr) \,\mathrm{d}\varphi + \mathrm{o}_{\Prob}(1)
\\
 &\weak&\int\gamma_{\theta
_B}(\varphi)\mathbb{G}_{\ltchat}^\angle(\varphi) \,\mathrm{d}\varphi=:
\Theta^{\mathrm{MD}},
\end{eqnarray*}
where $\ltcein(\varphi) = \ltc(\cos\varphi, \sin\varphi)$,
$\ltchatein= \ltchat(\cos\varphi, \sin\varphi)$, $\gamma_{\theta
_B} ( \varphi) =A_{\theta_B}^{-1} \delta_{\theta_B}^\angle(
\varphi)$,
$\delta_{\theta}^\angle( \varphi)=\partial_\theta\ltc(\cos
\varphi,\sin\varphi,\theta)$,
$ \mathbb{G}_{\ltchat}^\angle(\varphi)= \mathbb{G}_{\ltchat}(\cos
\varphi, \sin\varphi)$, and
\[
A_{\theta_B}:= \int\delta_{\theta_B}^\angle(\varphi)\delta
^{\angle}_{\theta_B}(\varphi)^T + \partial_\theta\delta_{\theta
_B}^\angle(\varphi) \bigl(\ltcein(\varphi;\theta_B) -\ltcein(\varphi
)\bigr) \,\mathrm{d}\varphi,
\]
with $\ltcein(\varphi; \theta) = \ltc(\cos\varphi, \sin\varphi;
\theta)$.
The limiting variable $\Theta^{\mathrm{MD}}$ is centered normally distributed
with variance
\[
\sigma^2=\int_{[0,\uppi/2]^2} \gamma_{\theta_B}(\varphi)\gamma
_{\theta_B}(\varphi') r(\cos\varphi,\sin\varphi,\cos\varphi
',\sin\varphi') \,\mathrm{d}(\varphi,\varphi'),
\]
where $r$ denotes the covariance functional of the process $\mathbb
{G}_{\ltchat}$
defined in (\ref{gf}).
\end{theorem}

To make use of the latter result in statistical applications, we need
the quantiles of the limiting distribution.
We propose to use
the multiplier bootstrap discussed in the previous section. The following
theorem shows that the $\mathit{pdm}$ and $\mathit{dm}$ bootstraps yield a valid
approximation of the distribution of the random variable $\Theta^{\mathrm{MD}}$.
\begin{theorem}\label{theo:thetamdasym}
If the assumptions of Theorem~\ref{theo:thetamdasy}
hold and $\Gamma_n$ denotes either the process
$\alpha_n^{\mathit{pdm}}$ (Theorem~\ref{theo:alphanpdm}) or
$\alpha_n^{\mathit{dm}}$ (Theorem~\ref{theo:alphandm})
obtained by the $ \mathit{pdm}$- or $\mathit{dm}$-bootstrap, respectively, then
\[
\Theta^{\mathrm{MD},m}_n:= \int\gamma_{\thetamd} (\varphi) \Gamma
_n^{\angle}(\varphi) \,\mathrm{d}\varphi \weakP \Theta^{\mathrm{MD}},
\]
where $\Gamma_n^{\angle}(\varphi)=\Gamma_n(\cos\varphi,\sin
\varphi)$, $\gamma_{\thetamd}=\hat A_{\thetamd}^{-1} \delta
_{\thetamd}^\angle( \varphi)$, and
\[
\hat A_{\thetamd}:= \int\delta_{\thetamd}^\angle(\varphi)\delta
^{\angle}_{\thetamd}(\varphi)^T + \partial_\theta\delta_{\thetamd
}^\angle(\varphi) \bigl(\ltcein(\varphi;\thetamd) -\ltchatein(\varphi
)\bigr) \,\mathrm{d}\varphi.
\]
\end{theorem}
\begin{table}
\tablewidth=\textwidth
\tabcolsep=0pt
\caption{\textit{Simulated coverage probability of the confidence intervals
based on the $\mathit{pdm}$ bootstrap for $n=1000$. In the last two columns, the
parameter $k$ is chosen as 50 for the estimation of $\theta$, whereas
$k$ is chosen as 200 for the bootstrap approximation of
$\Gb_{\ltchat}$}}\label{tab:confint}
\begin{tabular*}{\textwidth}{@{\extracolsep{\fill}}lllllll@{}}
\hline
&\multicolumn{6}{l}{$k$}
\\[-5pt]
&\multicolumn{6}{l@{}}{\hrulefill}
\\
 & \multicolumn{2}{l}{50} & \multicolumn{2}{l}{200} &
\multicolumn{2}{l}{50/200}
\\[-5pt]
& \multicolumn{2}{l}{\hrulefill} & \multicolumn{2}{l}{\hrulefill} &
\multicolumn{2}{l@{}}{\hrulefill}
\\
$\lambda_{L}$ & 90\% & 95\% & 90\% & 95\% & 90\% & 95\%
\\ \hline
0.25 & 0.895 & 0.955 & 0.014 & 0.044 & 0.830 & 0.915
\\
0.5 & 0.893 & 0.936 & 0.779 & 0.882 & 0.888 & 0.934
\\
0.75 & 0.838 & 0.887 & 0.900 & 0.949 & 0.863 & 0.894
\\ \hline
\end{tabular*}
\end{table}

Based on this result, it is possible to construct asymptotic confidence regions
for the  parameter~$\theta$, as well as to test point
hypotheses regarding the
parameter. Table~\ref{tab:confint} presents a small simulation study
regarding the
finite-sample coverage probabilities of some confidence intervals for
the parameter of a
Clayton tail copula. This interval is defined as $KI_{1-\alpha
}=[\thetamd-k^{-1/2} \hat q_{1-\alpha/2}, \thetamd-k^{-1/2} \hat
q_{\alpha/2}]$, where $\hat q_\beta$ denotes\vspace*{1pt} the estimated $\beta
$-quantile of the distribution of $\Theta^{\mathrm{MD}}_n$ based on the
bootstrap approximation provided by
Theorem~\ref{theo:thetamdasym}. The sample size is $n=1000$, and $B=500$
bootstrap replications are used for calculating the quantiles.\vadjust{\goodbreak} All
coverage probabilities are calculated by
$1000$ simulation runs. The parameter of the Clayton tail copula
is chosen such that the tail dependence coefficient varies in the set
$\{1/4,2/4,3/4\}$.

To investigate the impact of the choice of $k$, we chose three
different scenarios: $k=50$, $k=200$, and two different values of $k$,
namely $k=50$ for the estimator $\thetamd$ and $k=200$ for the
bootstrap estimator of the quantiles $\hat q_\beta$. This choice was
motivated by the findings in Section~\ref{subsec:fstail},
which indicate that a smaller value of $k$ should be used in the
estimator~$\hat\Lambda_L$.

The tables reveal that there is no unique ``optimal'' choice for $k$.
For $\lambda_L=0.25$,
the best results are obtained for the scenario with $k=50$, followed by
the scenario with two different values of~$k$.
Compare these findings with the results of Section~\ref{subsec:fstail}.
For $k=200$, the large bias of $\thetamd$ (compare the left side of
Figure~\ref{pic:cov}) demonstrates that the true parameter does not
lie in the estimated confidence interval for more than $95\%$ of the
repetitions. For stronger tail dependence, $\lambda_L=0.5, 0.75$, the
choice $k=200$ yields better results, with almost
perfect coverage probabilities for $\lambda_L=0.75$. Also note that
the case with $k=50$ in the estimator $\hat\Lambda_L$ and $k=200$
in the corresponding bootstrap statistic does not yield improvements
with respect to the approximation of the coverage probability compared
with the case where $k=50$.

\begin{table}[b]
\tablewidth=\textwidth
\tabcolsep=0pt
\caption{\textit{Relative efficiency of the method-of-moments type estimate to
the minimum distance estimate}}\label{tab:releff}
\begin{tabular*}{\textwidth}{@{\extracolsep{\fill}}llllllllll@{}}
\hline
 $\lambda$ & 0.1 & 0.2 & 0.3 & 0.4 & 0.5 & 0.6 & 0.7 & 0.8 & 0.9
\\
$E_{\lambda}$ &0.86 & 0.87 & 0.89 & 0.92 & 0.97 & 1.03 & 1.14 & 1.38 &
2.11
\\
\hline
\end{tabular*}
\end{table}

%
%
\begin{rem}
As pointed out at the beginning of this section, there are two
alternative estimators for parametric classes of tail copulas.
De~Haan \textit{et~al.} \cite{dehnevpen2008} proposed a censored maximum likelihood estimator and
proved weak convergence to a normal distribution, which involves the
partial derivatives up to the sixth order of the stable tail dependence
function. Einmahl \textit{et~al.} \cite{einkraseg2008} proposed a method-of-moments type
estimator and proved a similar statement as given in Theorem~\ref
{theo:thetamdasy} for the minimum distance estimate. Table~\ref
{tab:releff} compares the asymptotic variances of the method-of-moments
and minimum distance estimators for the parameter $\theta$ in the
Clayton family chosen such
that the coefficient of tail dependence $\lambda$ varies in the set $\{
0.1,\ldots,0.9\}$. The calculated values 
$E_{\lambda}$ are defined as
\[
E_{\lambda}=\frac{\mbox{Asymptotic variance of the minimum distance
estimate}}{\mbox{Asymptotic variance of the moment type estimate}}.
\]
(Note that we were not able to obtain the asymptotic variances for the
censored maximum likelihood estimator, because of the complicated
structure of the limiting distribution.) The method-of-moments
estimator requires the specification of a function $g$, which was
chosen as in Einmahl \textit{et~al.} \cite{einkraseg2008} as the indicator of the set $\{
\vect{x}\in[0,1]^2\dvt x_1+x_2\leq1\}$. We observe that none of the two
estimates is globally preferable over the other. For small amounts of
tail dependence, the minimum distance estimate performs slightly
better, whereas for increasing tail dependence, the method-of-moments
type estimator is more qualified from an asymptotic standpoint.

Also of note, the $\mathit{dm}$ and $\mathit{pdm}$ bootstraps can be used to construct a
consistent
approximation of the asymptotic distribution of the censored likelihood and
moment estimator investigated in de~Haan \textit{et~al.} \cite{dehnevpen2008} and
Einmahl \textit{et~al.} \cite{einkraseg2008}.
The main argument for proving consistency is that the limiting
distribution of the method-of-moments and minimum distance estimators
can be represented in the form $\Phi(\Gb_{\ltchat}, \ltc, \partial
\ltc)$ for some appropriate functional $\Phi$ depending on the method
of estimation. Here $\partial\ltc$ denotes any vector of partial
derivatives of $\ltc$ with respect to its coordinates or the
parameter. Given that the functional $\Phi$ is suitable smooth, the
bootstrap approximation is then obtained by $\Phi(\alpha_n, \ltchat,
\widehat{\partial\ltc})$, where $\alpha_n$ is $\alpha_n^{\mathit{pdm}}$ of
$\alpha_n^{\mathit{dm}}$ and $\widehat{\partial\ltc}$ is a consistent
estimate of $\partial\ltc$.%
\end{rem}

In what follows, we use the multiplier bootstrap to construct a computationally
efficient goodness-of-fit test for the hypothesis that the lower tail
copula has a
specific parametric form, that is,
%
\begin{equation} \label{h0class}
\mathcal{H}_0\dvt \ltc\in\mathcal{L}=\{ \ltc(\cdot, \theta) \vert
\theta\in\Theta\}, \qquad\mathcal{H}_1\dvt \ltc\notin\mathcal{L}.
\end{equation}
This problem was also been discussed by de~Haan \textit{et~al.} \cite{dehnevpen2008} and
Einmahl \textit{et~al.} \cite{einkraseg2008},
who proposed a comparison between a nonparametric estimate and a
parametric estimate of the
lower tail copula by an $L^2$-distance. In both cases, the limiting
distribution of the
corresponding test statistic under the null hypothesis depends in a
complicated way
on the process $\Gb_{\ltchat}$ and the unknown true parameter $\theta
_B$. Although
Einmahl \textit{et~al.} \cite{einkraseg2008} did not propose any bootstrap approximation,
de~Haan \textit{et~al.} \cite{dehnevpen2008}
proposed using the parametric bootstrap. However, Kojadinovic and Yan \cite{kojayan2010}
and Kojadinovic \textit{et~al.} \cite{kojyanhol2010} pointed out that for copula models,
approximations based on multiplier
bootstraps are computationally more efficient, especially for large
sample sizes.
We now illustrate how the multiplier bootstrap can be successfully
applied to the
problem of testing the hypothesis (\ref{h0class}).

To be precise, we propose comparing a parametric estimate\vspace*{1pt} (using the
minimum distance
estimate $\hat\theta_n^{\mathrm{MD}}$) and a nonparametric estimate of the
tail copula,
and rejecting the null hypothesis (\ref{h0class}) for large values of
the statistic
\[
\mathit{GOF}_n:= k \varrho(\ltchat, \ltc(\cdot;\thetamd))=k \int \bigl(
\ltchatein(\varphi)- \ltcein(\varphi;\thetamd) \bigr)^2 \,\mathrm{d}\varphi,
\]
where $\thetamd$ denotes the minimum distance estimate. If the
standard assumptions
of minimum distance estimation are satisfied (see page 89 in
B{\"{u}}cher \cite{buecher2011} for details), then we obtain for the process
$H_n = \sqrt{k} ( \ltchat-\ltc(\cdot;\thetamd) )$ under the null hypothesis
$\mathcal{H}_0\dvt \ltc=\ltc(\cdot;\theta_B)$
\begin{eqnarray*}
H_n  &=& \sqrt{k} \bigl( \ltchat- \ltc - \delta_\theta( \thetamd-\theta
) \bigr) + \mathrm{o}_{\Prob}(1)
\\
&=& \sqrt{k} \biggl( \ltchat- \ltc - \delta_\theta\int\gamma_\theta
(\varphi) \bigl(\ltchatein(\varphi) - \ltcein(\varphi)\bigr)\,\mathrm{d}\varphi\biggr)
+\mathrm{o}_{\Prob}(1)
\\
&\weak&\mathbb{G}_{\ltchat}-\delta_\theta\int\gamma_\theta
(\varphi) \mathbb{G}_{\ltchat}^\angle(\varphi) \,\mathrm{d}\varphi= \mathbb
{G}_{\ltchat} -
\delta_\theta\Theta^{\mathrm{MD}}.
\end{eqnarray*}
Under the alternative hypothesis, we get an additional summand,
\[
H_n = \sqrt{k} \bigl( \ltchat- \ltc - \delta_\theta( \thetamd-\theta) -
\bigl(\ltc(\cdot;\theta_B) -\ltc\bigr) \bigr) + \mathrm{o}_{\Prob}(1) ,
\]
which converges to either plus or minus infinity whenever
$\Lambda_L (\vect{x}, \theta_B) \neq\Lambda_L(\vect{x})$.
The continuous mapping theorem yields the following result.
\begin{theorem}\label{theo:gofasy}
Assume that assumptions of Theorem~\ref{theo:thetamdasy} are
satisfied. If the null hypothesis is valid, then
%
\begin{equation} \label{zproc}
\mathit{GOF}_n = \int\{ H_n^\angle(\varphi)\}^2 \,\mathrm{d}\varphi\weak Z:= \int
\bigl(\mathbb{G}_{\ltchat}^\angle(\varphi)-\delta_\theta^\angle
(\varphi) \Theta^{\mathrm{MD}} \bigr)^2 \,\mathrm{d}\varphi,
\end{equation}
whereas under the alternative,
$\mathit{GOF}_n = \int\{ H_n^\angle(\varphi)\}^2 \,\mathrm{d}\varphi\Pconv\infty$.
\end{theorem}

The critical values of the test, which rejects the null hypothesis for
large values of $\mathit{GOF}_n$,
can be calculated based on the following theorem. For a proof, see
Theorem 4.10 in B{\"{u}}cher \cite{buecher2011}.
\begin{theorem}\label{theo:gofboot}
If the assumptions of Theorem~\ref{theo:thetamdasy} hold and $\Gamma
_n$ denotes
either the process $\alpha_n^{\mathit{pdm}}$ (Theorem~\ref{theo:alphanpdm}) or
$\alpha_n^{\mathit{dm}}$ (Theorem~\ref{theo:alphandm})
obtained by the $ \mathit{pdm}$ bootstrap and the $\mathit{dm}$ bootstrap, respectively, then,
independent of the hypotheses, it holds that
\[
H_n^m:= \Gamma_n-\delta_{\thetamd} \int\gamma_{\thetamd}(\varphi
) \Gamma_n^{\angle}(\varphi) \,\mathrm{d}\varphi \weakP \mathbb{G}_{\ltchat
}-\delta_{\theta_B} \Theta^{\mathrm{MD}}.
\]
Therefore,
$\mathit{GOF}_n^m = \int\{ H_n^{m\angle}(\varphi)\}^2 \,\mathrm{d}\varphi \weakP Z$,
where $Z$ is as defined in \eqref{zproc}.
\end{theorem}
%
\begin{table}
\tablewidth=\textwidth
\tabcolsep=0pt
\caption{\textit{Simulated rejection probabilities of the $\mathit{pdm}$ bootstrap test
\protect\hyperref[testpar]{\emph{(4.9)}}
for the hypothesis \protect\hyperref[h0class]{\emph{(4.7)}}. The first three rows represent
models from the null hypothesis, and the next eight rows represent
alternative models. The sample size is $n=1000$, and $B=500$ bootstrap
replications were performed. $\lambda_L$ denotes the lower tail
dependence coefficient}}\label{tab4}
\begin{tabular*}{\textwidth}{@{\extracolsep{\fill}}lllllll@{}}
\hline
& \multicolumn{6}{l}{$k$}
\\[-5pt]
& \multicolumn{6}{l@{}}{\hrulefill}
\\
 & \multicolumn{3}{l}{50} & \multicolumn
{3}{l}{200}\\[-5pt]
& \multicolumn{3}{l}{\hrulefill} & \multicolumn
{3}{l@{}}{\hrulefill}
\\
\multicolumn{1}{l}{Model}& $\alpha=0.15$ & $\alpha=0.1$ & $\alpha
=0.05$ & $\alpha=0.15$ & $\alpha=0.1$ & $\alpha=0.05$
\\ \hline
Clayton ($\lambda_L=0.25$) & 0.124 & 0.087 & 0.037 & 0.174 & 0.108 &
0.049
\\
Clayton ($\lambda_L=0.5$) & 0.097 & 0.068 & 0.032 & 0.117 & 0.073 &
0.039
\\
Clayton ($\lambda_L=0.75$) & 0.091 & 0.048 & 0.018 & 0.091 & 0.058 &
0.024
\\[3pt]
Convex ($\lambda_L=1/12$) & 0.095 & 0.052 & 0.017 & 0.386 & 0.291 &
0.179
\\
Convex ($\lambda_L=2/12$) & 0.124 & 0.066 & 0.029 & 0.502 & 0.401 &
0.253
\\
Convex ($\lambda_L=3/12$) & 0.298 & 0.200 & 0.088 & 0.880 & 0.828 &
0.700
\\
Aneglog ($\lambda_L=0.2$) & 0.119 & 0.071 & 0.028 & 0.257 & 0.185 &
0.109
\\
Aneglog ($\lambda_L=0.4$) & 0.241 & 0.174 & 0.105 & 0.625 & 0.534 &
0.416
\\
Aneglog ($\lambda_L=0.6$) & 0.874 & 0.833 & 0.732 & 1.000 & 1.000 &
1.000
\\
Mixed ($\lambda_L=0.1$) & 0.118 & 0.069 & 0.022 & 0.523 & 0.424 & 0.268
\\
Mixed ($\lambda_L=0.3$) & 0.148 & 0.068 & 0.032 & 0.405 & 0.315 & 0.187
\\ \hline
\end{tabular*}
\end{table}

To investigate the finite-sample properties of a goodness-of-fit test
on the basis of the multiplier bootstrap, Table~\ref{tab4} presents
the simulated rejection probabilities of the $\mathit{pdm}$ bootstrap test,
%
\begin{equation} \label{testpar}
\mathit{GOF}_n > q^{(\mathit{pdm})}_{1- \alpha},
\end{equation}
where $q^{(\mathit{pdm})}_{1- \alpha}$ denotes the $(1- \alpha)$ quantile of
the bootstrap distribution. For the
null hypothesis, we considered the family of Clayton tail copulas as
the parametric class. In particular, we investigated three scenarios
corresponding to
a coefficient of tail dependence $\ltc(1,1)$ varying in $\{0.25, 0.5,
0.75\}$.
Under the alternative, we consider three models:
\begin{enumerate}[(3)]
\item[(1)] A convex
combination of the independence copula $\Pi(\vect{u})=u_1u_2$ and a
Clayton copula (with convex
parameter $1/3$), such that the tail copula is given by $\ltc(\vect
{x})=1/3 (x_1^{-\theta}+x_2^{-\theta})^{-1/\theta}$.
The parameter $\theta$ is chosen such that $\lambda_L=\ltc
(1,1)=1/3\times2^{-1/\theta}$ varies in the set $\{1/12,2/12,3/12\}$.
\item[(2)] The asymmetric negative logistic model (see
Joe \cite{joe1990}), defined by
\[
\ltc(1-t,t)= \bigl\{\bigl(\psi_1(1-t)\bigr)^{-\theta} + (\psi_2t)^{-\theta} \bigr\}
^{-1/\theta},\qquad t\in[0,1],
\]
with parameters $\psi_1=2/3, \psi_2=1$ and $\theta\in(0,\infty)$
chosen such that $\lambda_L = \ltc(1,1)$ varies in the set $\{
0.2,0.4,0.6\}$.
\item[(3)] The mixed model (see Tawn \cite{tawn1988}), given by
\[
\ltc(1-t,t)=\theta t(1-t), \qquad t\in[0,1],
\]
where the parameter $\theta\in[0,1]$ is chosen such that $\lambda_L
= \ltc(1,1)=\theta/2$ equals $0.1$ or $0.3$.
\end{enumerate}

The results are based on $1000$ simulation runs,
a sample size of $n=1000$, and two choices, $k=50, 200$, of the
parameter $k$. For each scenario, the critical values were calculated by
$B=500$ bootstrap replications with $\mathcal{U}(\{0,2\}
)$-multipliers. We observe a reasonable power and approximation of the
nominal level
in most cases. Under the null hypothesis, the test is conservative,
and this effect is increasing with the level of tail dependence.
For the mixed model with $k=50$, the power of the test is close to the
nominal level. This observation can be explained
by the fact that for $\lambda_L=0.5$ (which corresponds to the case
where $\theta=1$), the model is exactly the same as the Clayton
model with parameter $1$; that is, we get close to the null hypothesis
with increasing tail dependence.
Finally, we note that a choice of larger $k$ results in substantially
better power properties, whereas we noted no notable differences in the
quality of the approximation of the nominal level. Again, this may be
explained by the fact\vspace*{1pt} that bias terms in $\mathit{GOF}_n$ cancel out when
calculating the difference $ H_n = \sqrt{k} ( \ltchat-\ltc(\cdot
;\thetamd) )$. Therefore, we propose using rather large values for $k$
in applications of the goodness-of-fit test.
%
%
%
%
%

%
%

\begin{appendix}\label{append}
\section{Proofs}\label{appa}
\subsection{\texorpdfstring{Proof of Theorem \protect\ref{theo:ltchat}}{Proof of Theorem 2.2}}
Let $\Bc_\infty(\R_+)$ denote the set of functions $f\dvtx\R
_+\rightarrow\R$ (where $\R_+=[0,\infty)$) that are uniformly
bounded on compact sets (equipped with the topology of uniform
convergence on compact sets), and define $\mathcal{B}_\infty^{I,0}(\R
_+)$ as the subset of all nondecreasing functions $f\dvtx \R_+ \rightarrow
\R_+ $ that satisfy $f(0+)=0$ and for which $\sup\ran f<\infty$
(with $\ran f$ denoting the range of $f$) implies that there exists a
$x_0$ with $f(x_0)=\sup\ran f$. The latter condition implies that the
adjusted generalized inverse function, defined by
\[
f^-(z)=
\cases{
\sup\{x\in\R_+ | f(x)=0\},&\quad $z=0$,
\cr
\inf\{x\in\R_+ | f(x)\geq z \}, & \quad $0 < z < \sup\ran f$,
\cr
\inf\{x\in\R_+ | f(x)=\sup\ran f\}, &\quad $z\geq\sup\ran f$,}%
\]
remains in $\Bc_\infty(\R_+)$ for every $f\in\mathcal{B}_\infty
^{I,0}(\R_+)$.
Further, set
\[
\mathcal{B}_\infty^{I,0} (\bar{\R}^2_+) := \{ \gamma\in\mathcal
{B}_\infty(\bar{\R}^2_+) | \gamma(\cdot,\infty)
\in\mathcal{B}_\infty^{I,0}(\R_+), \gamma(\infty,\cdot)\in
\mathcal{B}_\infty^{I,0}(\R_+)\}
\]
and now define a map $\Phi\dvtx \mathcal{B}_\infty^{I,0} (\bar{\R
}^2_+) \rightarrow \mathcal{B}_\infty(\bar{\R}^2_+)$ by
\[
\gamma\longmapsto \Phi(\gamma)(x,y)= \cases{
 \gamma( \gamma^-(x,\infty),
\gamma^-(\infty,y) ), &\quad if $x,y\ne\infty$,
\cr
 \gamma(\gamma^-(x,\infty),
\infty), &\quad if $y=\infty$,
\cr
 \gamma( \infty, \gamma^-(\infty,y) ), &\quad if
$x=\infty$}
\]
(see also Schmidt and Stadtm{\"{u}}ller \cite{schmstad2006}). Observing that $\ltctilde\in
\mathcal{B}_\infty^{I,0} (\bar{\R}^2_+)$, and that
the adjusted generalized inverse of $\ltctilde(x,\infty)$ is
given by $\frac{n}{k}F_1(F_{n1}^-(kx/n))$, we can
conclude that $\Phi(\ltc)=\ltc$ and $\Phi(\ltctilde) = \ltchat$
($\Prob$-almost surely), and the proof of Theorem~\ref{theo:ltchat}
follows from the functional delta method (Theorem 3.9.4 of {V}an~der Vaart and Wellner~\cite{vandwell1996}) and the following Lemma, which is an extension of the
result in the proof of Theorem 5 of Schmidt and Stadtm{\"{u}}ller \cite{schmstad2006} to our weaker
conditions.

\begin{lemma}\label{lem:hadamard}
Let $\ltc$ be a lower tail copula whose partial derivatives satisfy
the first-order property \eqref{fotc} for $p=1,2$.
Then $\Phi$ is Hadamard-differentiable at $\ltc$ tangentially to the set
\[
\mathcal{C}^0 (\bar{\R}^2_+) = \{ \gamma\in\mathcal{B}_\infty
(\bar{\R}^2_+) | \gamma\mbox{ continuous with } \gamma(\cdot,0 )=
\gamma(0,\cdot)=0 \}.
\]
Its derivative at $\ltc$ in $\gamma\in\mathcal{C}^0 (\bar{\R
}^2_+) $ is given by
%
\begin{equation} \label{Phi'}
\Phi'_{\ltc}(\gamma)(\vect{x})=\gamma(\vect{x}) - \partial_1\ltc
(\vect{x})\gamma(x_1,\infty)-\partial_2\ltc(\vect{x})\gamma
(\infty,x_2),
\end{equation}
where $\partial_p \ltc, p=1,2$ is defined as 0 on the set $\{\vect
{x}\in\bar\R_+^2 | x_p\in\{0,\infty\}\}$.
\end{lemma}

\begin{pf}
Decompose $ \Phi= \Phi_3\circ\Phi_2 \circ\Phi_1 $,
where
\begin{eqnarray*}
&&\Phi_1\dvtx   \mathcal{B}_\infty^{I,0} (\bar{\R}^2_+) \rightarrow
\mathcal{B}_\infty^{I,0} (\bar{\R}^2_+)\times\mathcal{B}_\infty
^{I,0} ({\R}_+)\times\mathcal{B}_\infty^{I,0} ({\R}_+),
\\
&&\quad \gamma\longmapsto(\gamma, \gamma(\cdot,\infty), \gamma(\infty
,\cdot) ),
\\
&&\Phi_2 \dvtx \mathcal{B}_\infty^{I,0} (\bar{\R}^2_+)\times\mathcal
{B}_\infty^{I,0} ({\R}_+)\times\mathcal{B}_\infty^{I,0} ({\R}_+)
\rightarrow\mathcal{B}_\infty^{I,0} (\bar{\R}^2_+)\times\mathcal
{B}_\infty^{I,0,-} ({\R}_+)\times\mathcal{B}_\infty^{I,0,-} ({\R
}_+),
\\
&&\quad (\gamma, f,g ) \longmapsto(\gamma, f^-, g^-) ,
\\
&&\Phi_3\dvtx
\mathcal{B}_\infty^{I,0} (\bar{\R}^2_+)\times\mathcal{B}_\infty
^{I,0,-} ({\R}_+)\times\mathcal{B}_\infty^{I,0,-} ({\R}_+)
\rightarrow\mathcal{B}_\infty(\bar{\R}^2_+),
\\
&&\quad (\gamma,f,g ) \longmapsto \cases{
\gamma( f(x), g(y) ), &\quad if $x,y\ne
\infty$,
\cr
 \gamma(f(x), \infty), & \quad if $y=\infty$,
 \cr
  \gamma( \infty, g(y) ), & \quad if
 $x=\infty$,}
\end{eqnarray*}
where $\mathcal{B}_\infty^{I,0,-} ({\R}_+) $ denotes the set of all
adjusted generalized inverse functions $f^-$ with $f\in\mathcal
{B}_\infty^{I,0} ({\R}_+)$.
Now $\Phi_1$ is Hadamard-differentiable at $\ltc$ tangentially to
$\mathcal{C}^0 (\bar{\R}^2_+)$, because it is linear and continuous.
The second map, $\Phi_2$, is Hadamard-differentiable at $(\ltc,\id
_{\R_+}, \id_{\R_+})$ tangentially to
$\mathcal{C}^0 (\bar{\R}^2_+) \times\mathcal{C}^0 (\R_+) \times
\mathcal{C}^0 (\R_+) $,
where $\mathcal{C}^0 (\R_+) $ consists of all continuous functions
$f$ on $\R_+$ with $f(0)=0$ and
its derivative at $(\ltc,\id_{\R_+},\id_{\R_+})$ in $(\gamma, f,
g)$ is given
by $\Phi'_{2, (\ltc,\id_{\R_+},\id_{\R_+}) }(\gamma,f,g) =
(\gamma, -f,-g)$. The proof follows along similar lines as the proof
of Theorem 5 in Schmidt and Stadtm{\"{u}}ller \cite{schmstad2006}, page 321, and thus is omitted; we
simply note that $(\id_{\R_+}+t_nf_n)^-(x)>0$ for all $x>0$ is
implied by the additional assumption of continuity in $0$ for functions
in the set $\mathcal{B}^{I,0}(\R_+)$.
More effort is needed to show
that $\Phi_3$ is Hadamard-differentiable at $(\ltc,\id_{\R_+}, \id
_{\R_+})$ tangentially to $\mathcal{C}^0 (\bar{\R}^2_+) \times
\mathcal{C}^0 (\R_+) \times\mathcal{C}^0 (\R_+)$ with derivative
\[
\Phi'_{3, (\ltc,\id_{\R_+},\id_{\R_+}) }(\gamma,f,g)(\vect{x})
= \gamma(\vect{x}) + \partial_1 \ltc(\vect{x}) f(x_1) + \partial
_2 \ltc(\vect{x})g(x_2).
\]
To see this, let $t_n\rightarrow0$, $(\gamma_n, f_n, g_n)\in\mathcal
{B}_\infty(\bar{\R}^2_+)\times\mathcal{B}_\infty({\R}_+)\times
\mathcal{B}_\infty({\R}_+)$
with $(\gamma_n,f_n,g_n)\rightarrow(\gamma,f,g) \in\mathcal{C}^0
(\bar{\R}^2_+) \times\mathcal{C}^0 (\R_+) \times\mathcal{C}^0
(\R_+) $
such that $(\ltc+t_n\gamma_n, \id_{\R_+}+t_nf_n, \id_{\R_+} +
t_ng_n)\in\mathcal{B}_\infty^{I,0} (\bar{\R}^2_+)\times\mathcal
{B}_\infty^{I,0,-} ({\R}_+)\times\mathcal{B}_\infty^{I,0,-} ({\R}_+)$.
Now $\Phi_3$ is linear in its first argument, and we introduce the
decomposition
\[
t_n^{-1} \{ \Phi_3(\ltc+t_n\gamma_n, \id_{\R_+}+t_nf_n, \id_{\R
_+} + t_ng_n) - \Phi_3(\ltc,\id_{\R_+},\id_{\R_+}) \} = L_{n1}+L_{n2},
\]
where
\begin{eqnarray*}
L_{n1} &=& t_n^{-1} \{ \Phi_3(\ltc, \id_{\R_+}+t_nf_n, \id_{\R_+}
+ t_ng_n) - \Phi_3(\ltc,\id_{\R_+},\id_{\R_+}) \},
\\
L_{n2} &=& \Phi_3(\gamma_n, \id_{\R_+}+t_nf_n, \id_{\R_+} + t_ng_n).
\end{eqnarray*}
By the definition of $d$, it suffices to show uniform convergence on
the  sets $T_i$, $i\in\Nat$. Because $T_i\subset\bar\R
_+^2$ is compact, $(f_n,g_n)$ converges uniformly and $\gamma$ is
uniformly continuous; thus $L_{n2}$ uniformly converges to $\gamma$.

Considering $L_{n1}$, we split the investigation into six different
cases, depending on the position of $\vect{x}\in T_i$. First, let
$\vect{x}\in(0,i]^2$. A series expansion at $\vect{x}$ yields
\[
L_{n1}=\partial_1\ltc(\vect{x})f_n(x_1) + \partial_2\ltc(\vect
{x}) g_n(x_2) + r_n(\vect{x}),
\]
where the error term $r_n$ can be written as
\[
r_n(\vect{x}) = \bigl(\partial_1\ltc(\vect{y})- \partial_1\ltc(\vect
{x}) \bigr) f_n(x_1) + \bigl(\partial_2\ltc(\vect{y})-\partial_2\ltc(\vect
{x}) \bigr)g_n(x_2)
\]
with some intermediate point $\vect{y}=\vect{y}(n)$ between $\vect
{x}$ and $(x_1+t_nf_n(x_1), x_2+t_nf_n(x_2))$.
The dominating term converges uniformly to $\partial_1 \ltc(\vect
{x}) f(x_1) + \partial_2 \ltc(\vect{x})g(x_2)$; thus it remains to
show that $r_n(\textbf{x})$ converges to $0$ uniformly in $\vect{x}$.
For a given $\varepsilon>0$,
uniform convergence of $f_n$ and uniform continuity of $f$ on $[0,i]$,
as well as the fact that $f(0)=0$, allows us to choose a $\delta>0$
such that $|f_n(x_1)|<\varepsilon$ for all $x_1<\delta$. Because
partial derivatives of tail copulas are bounded by $1$, the first
term of $r_n(\textbf{x})$ is uniformly small for $x_1<\delta$. On the
quadrangle $[\delta,i]\times(0,i]$, the partial derivative $\partial
_1\ltc$ is
uniformly continuous, which yields the desired convergence under
consideration of $\vect{y}(n)\rightarrow\vect{x}$ and boundedness of $f$.
The same arguments apply for the second derivative, and the case $\vect
{x}\in(0,i]^2$ is finished.

Now consider the case $\vect{x}\in(0,i]\times\{0\}$. By Lipschitz
continuity of $\ltc$ on $\R_+^2$ (see Theorem 1 in
Schmidt and Stadtm{\"{u}}ller \cite{schmstad2006}), we get
\begin{eqnarray*}
|L_{n1}(x_1,0)| &=& t_n^{-1} \bigl|\ltc\bigl(x_1+t_nf_n(x_1),t_ng_n(0)\bigr)\bigr |
\\
&=& t_n^{-1}\bigl|\ltc\bigl(x_1+t_nf_n(x_1),t_ng_n(0)\bigr) - \ltc
\bigl(x_1+t_nf_n(x_1),0\bigr)\bigr|
\\
&\leq&|g_n(0)| \rightarrow g(0)=0.
\end{eqnarray*}
Because $\partial_1\ltc(x_1,0)f(x_1)+\partial_2\ltc(x_1,0)g(0)=0$,
this yields the assertion. The arguments are similar for the cases
$\vect{x}=(0,0)^T$ and $\vect{x}\in\{0\}\times(0,i]$, and we
proceed with $\vect{x}\in[0,i]\times\{\infty\}$ (and analogously
$\vect{x}\in\{\infty\}\times[0,i]$)
\[
L_{n1}(x_1,\infty) = t_n^{-1}\bigl(\ltc\bigl(x_1+t_nf_n(x_1),\infty\bigr) - \ltc
(x_1,\infty)\bigr) = f_n(x_1)\rightarrow f(x_1).
\]
This yields the assertion by $\partial_1\ltc(x_1,\infty)=1$ and
$\partial_2\ltc(x_1,\infty)=0$. To conclude, $\Phi_3$ is
Hadamard-differentiable as asserted.

An application of the chain rule (see Lemma 3.9.3 in
{V}an~der Vaart and Wellner \cite{vandwell1996}) completes the proof of the lemma.
\end{pf}
\subsection{\texorpdfstring{Proof of Theorem \protect\ref{theo:atildem}}{Proof of Theorem 3.1}} \label{subsecproof}
In light of Lemma~\ref{bcond} in Appendix \ref{app:B} (an analog of Theorem
1.6.1 in \cite{vandwell1996} for the case of conditional weak
convergence), the proof of conditional weak
convergence of $\tilde\alpha_n^m$ in $\mathcal{B}_\infty(\bar{\R}_+^2)$
can be given for each $\ell^\infty(T_i)$ separately. For brevity, we
suppress the index $i$ and write $T=T_i$.
Recalling the notation of $f_{n,\vect{x}} (\vect{U}_i)$ in \eqref{f_nx}
we can express $\tilde\alpha_n^m$ as
\[
\tilde{\alpha}^{m}_n(\vect{x})=\frac{\mu}{\tau}\sqrt{k} (
\ltctilde^\xi-\ltctilde)= \frac{\mu}{\tau}\frac{1}{\sqrt
{n}}\sum_{i=1}^n
\biggl( \frac{\xi_i}{\bar{\xi}_n}-1\biggr )f_{n,\vect{x}}(U_i),
\]
and the assertion now follows by an application of Theorem 11.23 in
Kosorok \cite{kosorok2008}.
For this purpose, we show that the assumptions for
this result are satisfied. Let $\mathcal{F}_n=\{ f_{n,\vect{x}} \dvt
\vect{x}\in T\}$
be a class of functions changing with $n$ and let
\[
F_n(\vect{u})=\sqrt{\frac{n}{k}} \ind \{ u_1 \leq k{ i}/n \mbox{
or } u_2 \leq ki/n \},
\]
%
denote a corresponding sequence of envelopes of
$\mathcal{F}_n$. We must prove the following:
\begin{enumerate}[(iii)]
\item[(i)] $(\mathcal{F}_n,F_n)$ satisfies the bounded
uniform entropy integral condition
%
\begin{equation} \label{BUIE}
\limsup_{n\rightarrow\infty} \sup_{Q}\int_0^1 \sqrt{\log N(\eps
\Vert F_n\Vert_{Q,2}, \mathcal{F}_n, L_2(Q))} \,\mathrm{d}\eps<\infty,
\end{equation}
where for each $n$, the supremum ranges over all probability measures $Q$
with finite support and $\Vert F_n\Vert_{Q,2}= (\int F_n(x)^2 \,\mathrm{d}Q(x) )^{1/2}>0$.
\item[(ii)] The limit $H(\vect{x},\vect{y})=\lim
_{n\rightarrow\infty} \Exp[ \tilde\alpha_n(\vect{x})\tilde\alpha
_n(\vect{y})]$
exists for every $\vect{x}$ and $\vect{y}$ in $T$.
\item[(iii)] $\limsup_{n\rightarrow\infty} \Exp
F_{n}^2(\vect{U}_1) < \infty$
\item[(iv)] $\lim_{n\rightarrow\infty} \Exp F_{n}^2(\vect
{U}_1) \ind\{ F_{n}(\vect{U}_1)>\eps\sqrt{n} \} = 0$ for all $\eps>0$.
\item[(v)] $\lim_{n\rightarrow\infty} \rho_n(\vect
{x},\vect{y})=\rho(\vect{x},\vect{y})$
for all $\vect{x},\vect{y}\in T$, where
%
\begin{equation}
\rho_n(\vect{x},\vect{y})= \bigl( \Exp\bigl(f_{n,\vect{x}}(U_1) -
f_{n,\vect{y}}(U_1) \bigr)^2 \bigr)^{1/2}.
\end{equation}
Furthermore, for all sequences $(\vect{x}_n)_n, (\vect{y}_n)_n$ in
$T$, the convergence
$\rho_n(\vect{x}_n,\vect{y}_n)\rightarrow0$ holds, provided $\rho
(\vect{x}_n,\vect{y}_n)\rightarrow0$.
\item[(vi)] The sequence $\mathcal{F}_n$ of classes is almost
measurable
Suslin (AMS), that is, for all $n\geq1$ there exists a Suslin
topological space
$T_n\subset T$ with Borel sets $\mathcal{B}_n$ such that
\begin{enumerate}[(b)]
\item[(a)] $\Prob^*(\sup_{\vect{x} \in T} \inf_{\vect{y}\in T_n}
|f_{n,\vect{x}}(\vect{U}_1)-f_{n,\vect{y}}(\vect{U}_1)| > 0 ) = 0$,
\item[(b)] $f_{n,\cdot}\dvtx[0,1]^2\times T_n\rightarrow\R$ is
$\mathcal{B}|_{[0,1]^2}\times\mathcal{B}_n$-measurable for
$i=1,\ldots,n$.
\end{enumerate}
\end{enumerate}
To prove the bounded uniform entropy integral condition (i), we
decompose $\mathcal{F}_n=\bigcup_{i=1}^3 \mathcal{F}_n^{(i)}$ with
$\mathcal{F}_n^{(i)}=\{ f_{n,\vect{x}}^{(i)}, \vect{x}\in T\} $ and
\begin{eqnarray*}
f_{n,\vect{x}}^{(1)}(\vect{U}_i) &=& \sqrt{\frac{n}{k}} \ind \{
U_{i1} \leq kx_1/n \} \ind\{ x_2=\infty\},
\\
f_{n,\vect{x}}^{(2)}(\vect{U}_i)  &=& \sqrt{\frac{n}{k}} \ind \{
U_{i2} \leq kx_2/n \} \ind\{ x_1=\infty\},
\\
f_{n,\vect{x}}^{(3)}(\vect{U}_i) &=& \sqrt{\frac{n}{k}} \ind \{
U_{i1} \leq kx_1/n, U_{i2} \leq kx_2/n \} \ind\{ x_1<\infty,
x_2<\infty\}.
\end{eqnarray*}
The corresponding envelopes of the classes $\mathcal{F}_n^{(i)}$ are
given by
\begin{eqnarray*}
F_n^{(1)}(\vect{U}_i) & =& \sqrt{\frac{n}{k}} \ind(U_{i1}\leq ki/n),\\
F_n^{(2)}(\vect{U}_i) &=& \sqrt{\frac{n}{k}} \ind(U_{i2}\leq ki/n),
\\
F_n^{(3)}(\vect{U}_i) & =& \sqrt{\frac{n}{k}} \ind(U_{i1}\leq ki/n,
U_{i2} \leq ki/n),
\end{eqnarray*}
so that $F_n(\vect{U}_i) = \max_{i=1}^3\{F_n^{(i)}(\vect{U}_i) \}$.
If we prove
that the sequences $(\mathcal{F}_n^{(i)},F_n^{(i)})$ satisfy the bounded
uniform integral entropy condition given in (\ref{BUIE}), then the condition
also holds for $(\mathcal{F}_n,F_n)$ by Lemma~\ref{BUIEunion} in the \hyperref[app:C]{Appendix},
and thus the assertion in (i) is proved. We consider only the (hardest)
case of
$\mathcal{F}_n^{(3)}$. Note that
${\mathcal{F}}^{(3)}_n=\{ {f}_{n,\vect{x}} , \vect{x}\in[0,i]^2\}
=\mathcal{G}_n^{(1)} \cdot\mathcal{G}_n^{(2)}$,
where
\begin{eqnarray*}
{f}_{n,\vect{x}} &=& (n/k)^{1/2} \ind\{ U_{i1}\leq kx_1/n, U_{i2}\leq
kx_2/n \},
\\
\mathcal{G}_n^{(j)} &=& \bigl\{ g_{n,t}=(n/k)^{1/4} \ind\{U_{ij}\leq kt/n
\} | t\in[0,i]\bigr\}
\end{eqnarray*}
for $j=1,2$. Because the functions $g_{n,t}$ are increasing in $t$, the
$G_n^{(j)}$ are
VC classes with VC index 2. Thus, by Lemma 11.21 in Kosorok \cite{kosorok2008},
both classes satisfy the bounded uniform integral entropy condition
(\ref{BUIE}).
Proposition 11.22 in Kosorok \cite{kosorok2008} shows that $\mathcal{F}_n^{(3)}$
has the same property, and by the discussion at the beginning of this paragraph,
(\emph{i}) is\vspace*{1pt} satisfied.

For the proof of (ii), note that
$\Exp[ \tilde\alpha_n(\vect{x}) \tilde\alpha_n(\vect{y})]=n/k (
C(\frac{(\vect{x}\wedge\vect{y})k}{n})-C(\frac{\vect{x}
k}{n})C(\frac{\vect{y} k}{n}) )$, which converges to $\ltc(\vect
{x}\wedge\vect{y})=:H(\vect{x},\vect{y})$, because $\frac
{n}{k}C(\frac{\vect{x} k}{n})C(\frac{\vect{y} k}{n})\rightarrow
0$.

Regarding (iii) and (iv), we note that $\Exp F_{n}(\vect{U}_1)^2=2i -
\frac{n}{k}C(ik/n,ik/n)$,
which converges to $2i-{\ltc}(i,i)$. Furthermore,
\begin{eqnarray*}
\Exp F_{n}^2(\vect{U}_1)\ind\bigl\{ F_{n}(\vect{U}_1) > \eps\sqrt{n} \bigr\}
&=& \int_{\{ F_{n}(\vect{U}_1)>\eps\sqrt{n} \} } F_{n}^2(\vect
{U}_1) \,\mathrm{d}\Prob
\\
&\leq&\frac{n}{k}\Prob\biggl(\frac{1}{{k}}\ind\{ U_{11}\leq k{ i}/n
\mbox{ or } U_{12} \leq ki/n\} >\eps\biggr) = 0
\end{eqnarray*}
for sufficiently large $n$, such that $k>1/\varepsilon$. For (v), we
note that
\begin{eqnarray*}
\rho_n(\vect{x},\vect{y}) &=& \bigl( \Exp\bigl(f_{n,\vect{x}}(\vect{U}_1) -
f_{n,\vect{y}}(\vect{U}_1) \bigr)^2 \bigr)^{1/2}
\\
&=&\sqrt{\frac{n}{k}} \bigl(C(\vect{x} k/n)-2C\bigl((\vect{x} \wedge\vect
{y})k/n\bigr)+C(\vect{y} k/n) \bigr)^{1/2}
\\
&\rightarrow& \bigl( {\ltc}(\vect{x})-2{\ltc}(\vect{x}\wedge\vect
{y})+{\ltc}(\vect{y}) \bigr)^{1/2} =: \rho(\vect{x},\vect{y}).
\end{eqnarray*}
In light of Theorem 1 in Schmidt and Stadtm{\"{u}}ller \cite{schmstad2006}, we have locally uniform
convergence
in the latter expression, which yields the second condition stated in
(v).

For the proof of condition (vi), we use Lemma 11.15 and the discussion on
page 224 in Kosorok \cite{kosorok2008} and show separability of $\mathcal
{F}_n$; that is, for every $n\geq1$, there exists a countable subset
$T_n\subset T$ such that
\[
\Prob^* \Bigl( \sup_{\vect{x}\in T} \inf_{\vect{y}\in T_n} |f_{n,\vect
{y}}(\vect{U}_1)-f_{n,\vect{x}}(\vect{U}_1)| >0 \Bigr)=0.
\]
Choose $T_n=(\Q\cap[0,i]\times\{\infty\}) \cup(\{\infty\} \times
\Q\cap[0,i]) \cup(\Q^2 \cap [0,i]^2)$. We then have (note that the
functions $f_{n,\vect{x}}$ are built by indicators)
that for every $\omega$ and every $\vect{x}\in T$, there is an $\vect
{y}\in T_n$
with $|f_{n,\vect{x}}(\vect{U}_1(\omega))-f_{n,\vect{y}}(\vect
{U}_1(\omega))|=0$.
This yields the assertion, and thus the proof of Theorem~\ref
{theo:atildem} is complete.

\begin{rem}\label{rem:ap}
Given that
\[
\tilde\alpha_n(\vect{x})=\frac{1}{\sqrt{n}}\sum_{i=1}^n
\bigl(f_{n,\vect{x}}(\vect{U}_i)-\Exp f_{n,\vect{x}}(\vect{U}_i)\bigr)
\]
and that in Section~\ref{subsecproof} we provided the sufficient
conditions for an application of Theorem 11.20 in Kosorok \cite{kosorok2008},
we obtain an alternative proof of Lemma~\ref{lem:ltctilde}.
\end{rem}
%
\subsection{\texorpdfstring{Proof of Theorem \protect\ref{theo:alphandm}}{Proof of Theorem 3.4}}
For technical reasons, we give a proof of Theorem~\ref{theo:alphandm}
in advance of the proofs of Theorems~\ref{theo:betan} and~\ref
{theo:alphanpdm}. The proof is essentially
a consequence of a bootstrap version of the functional delta method
(see Theorem 12.1
in Kosorok \cite{kosorok2008}).
Because this result holds only for Banach space-valued stochastic
processes, some adjustments must be made.
Note that the space $\mathcal{B}_\infty(\bar{\R}^2_+)$ is a
complete topological vector space with a metric $d$,
and some care is necessary whenever technical results depending on the
norm are used.

Given Lemma~\ref{lem:ltctilde} and
Theorem~\ref{theo:atildem}, we have
\[
\sqrt{k} ( \ltctilde-\ltc) \weak\Gb_{\ltctilde},\qquad
\sqrt{k} \frac{\mu}{\tau} (\ltctilde^\xi- \ltctilde ) \weakP\Gb
_{\ltctilde}
\]
in $\mathcal{B}_\infty(\bar{\R}^2_+) $.
Observing that
the generalized inverses\vspace{-1pt} of $\ltctilde(x,\infty)$ and $\ltctilde^\xi
(x,\infty)$ are ($\Prob$-almost surely)
given by $\frac{n}{k}F_1(F_{n1}^-(kx/n))$ and $\frac
{n}{k}F_1(F_{n1}^{\xi-}(kx/n))$, respectively, we can
conclude that $\Phi(\ltc)=\ltc, \Phi(\ltctilde) = \ltchat$ and
$\Phi(\ltctilde^\xi)=\ltchat^{\xi,\xi}$ ($\Prob$-almost surely).
By Lemma~\ref{lem:hadamard}, $\Phi$ is Hadamard-differentiable on
$\mathcal{B}_\infty^{I} (\bar{\R}^2_+)$ at $\gamma_0=\ltc$
tangentially to $\mathcal{C}^0 (\bar{\R}^2_+) \subset\mathcal
{B}_\infty(\bar{\R}^2_+)$. Therefore, it remains to argue why
Theorem 12.1
in Kosorok \cite{kosorok2008} can be applied in the present context.

A careful inspection of the proof of Theorem 12.1 in Kosorok
\cite{kosorok2008} shows that\vspace*{1pt} properties going beyond our specific
assumptions (i.e., the complete topological vector space $(\mathcal
{B}_\infty(\bar{\R}^2_+),d)$)
are used only three times.
First, the mapping $\Phi'_{\ltc}$ needs to be extended to the
whole\vspace*{-1pt}
space $\mathcal{B}_\infty(\bar{\R}^2_+)$, which is possible
by using equation (\ref{Phi'}) as the defining identity. Second, the
proof of Theorem~12.1 in Kosorok \cite{kosorok2008}
uses the usual functional delta method as stated in Theorem 2.8 in the
same reference, but this
result can be replaced by Theorem 3.9.4 in {V}an~der Vaart and Wellner \cite{vandwell1996}, which provides
a functional delta method that holds in general metrizable topological
vector spaces.
Finally, the proof of Theorem 12.1 in Kosorok \cite{kosorok2008}
makes use of a bootstrap continuous mapping theorem (see Theorem 10.8
in Kosorok \cite{kosorok2008}), which yields
\[
\sqrt{k}\frac{\mu}{\tau}(\ltctilde^\xi- \ltctilde) \weakP\Gb
_{\ltctilde}\quad \Rightarrow \quad\Phi'_{\ltc}\biggl(\sqrt{k}\frac{\mu}{\tau
}(\ltctilde^\xi- \ltctilde)\biggr) \weakP\Phi'_{\ltc}(\Gb_{\ltctilde}).
\]
In our specific context, this statement follows immediately from the
Lipschitz continuity of the derivative
$\Phi'_{\ltc}$ and an application of Lemma~\ref{lem:contboot} in
Appendix~\ref{app:C}.

\subsection{\texorpdfstring{Proof of Theorem \protect\ref{theo:betan}}{Proof of Theorem 3.2}}

Consider the mapping
$\Psi\dvtx \mathcal{B}_\infty^{I,0} (\bar{\R}^2_+)\times\mathcal
{B}_\infty^{I,0} (\bar{\R}^2_+) \longrightarrow\mathcal{B}_\infty
(\bar{\R}^2_+)$ defined by
$\Psi= \Phi_3\circ\Phi_2 \circ\Psi_1$,
where $\Phi_3$ and $\Phi_2$ are defined in the proof of Lemma~\ref
{lem:hadamard} and $\Psi_1$ is given by
\begin{eqnarray*}
&&\Psi_1\dvtx  \mathcal{B}_\infty^{I,0} (\bar{\R}^2_+)\times\mathcal
{B}_\infty^{I,0} (\bar{\R}^2_+) \rightarrow\mathcal{B}_\infty
^{I,0} (\bar{\R}^2_+)\times\mathcal{B}_\infty^{I,0} ({\R
}_+)\times\mathcal{B}_\infty^{I,0} ({\R}_+),
\\
&&\quad (\beta, \gamma) \longmapsto(\beta, \gamma(\cdot,\infty),
\gamma(\infty,\cdot) ).
\end{eqnarray*}
Note that we obtain $\Psi(\ltc,\ltc)=\ltc$, $\Psi(\ltctilde
,\ltctilde)=\ltchat$ and
$\Psi(\ltctilde^\xi,\ltctilde)=\ltchat^{\xi,\cdot}$ ($\Prob
$-almost surely).
Clearly, $\Psi_1$ is Hadamard-differentiable at $(\ltc,\ltc)$,
because it is linear and continuous. $\Phi_2$ and $\Phi_3$ are
Hadamard-differentiable tangentially to suitable subspaces as well (see
the proof of Lemma~\ref{lem:hadamard}). By an application of the chain
rule (see Lemma 3.9.3 in \cite{vandwell1996}), we can conclude that
$\Psi$ is Hadamard-differentiable
$(\ltc,\ltc)$ tangentially to $\mathcal{C}^0(\bar{\R}^2_+) \times
\mathcal{C}^0 (\bar{\R}^2_+)$ with derivative
%
\begin{equation}
\label{Psi'}
\Psi'_{(\ltc,\ltc)}(\beta,\gamma)(\vect{x})=\beta(\vect{x}) -
\partial_1\ltc(\vect{x})\gamma(x_1,\infty)-\partial_2\ltc(\vect
{x})\gamma(\infty,x_2).
\end{equation}

Note that, unlike in the previous proof, we do not have weak
convergence (resp., weak conditional convergence) of $\sqrt{k}
((\ltctilde,\ltctilde)-(\ltc,\ltc) )$ and
$ \frac{\mu}{\tau}\sqrt{k} ((\ltctilde^\xi,\ltctilde) -
(\ltctilde,\ltctilde) )$ toward the same limiting field, which is
necessary for an application
of the functional delta method to the bootstrap (see, e.g., Theorem
12.1 in Kosorok \cite{kosorok2008}). Nevertheless, we can mimic certain
steps in the proof of this theorem to
obtain the result. To be precise, note that, by arguments analogous to
those on page 236 of Kosorok \cite{kosorok2008}, we obtain that
\[
\sqrt{k}
\pmatrix{ \ltctilde^\xi-\ltc
\cr
 \ltctilde-\ltc}
\weak \pmatrix{
 c^{-1} \Gb_1 + \Gb_2
 \cr
  \Gb_2},
\]
unconditionally, where $\Gb_1$ and $\Gb_2$ denote independent copies
of $\Gb_{\ltctilde}$ and $c=\mu\tau^{-1}$.
Hadamard-differentiability of the mapping
$(\beta,\gamma)\mapsto(\Psi(\beta,\gamma),\Psi(\gamma,\gamma
),(\beta,\gamma),(\gamma,\gamma))$ and the usual functional delta method
(Theorem 3.9.4 in {V}an~der Vaart and Wellner \cite{vandwell1996}) yield
\[
\sqrt{k}
\pmatrix{
 \Psi(\ltctilde^\xi,\ltctilde) - \Psi(\ltc,\ltc)
\cr
 \Psi(\ltctilde,\ltctilde) - \Psi(\ltc,\ltc)\vspace*{2pt}
\cr
 (\ltctilde^\xi
,\ltctilde) -(\ltc,\ltc)
\cr (\ltctilde,\ltctilde) - (\ltc,\ltc)}
\weak
\pmatrix{
 \Psi'_{(\ltc,\ltc)}(c^{-1} \Gb_1 + \Gb_2 ,\Gb_2)\vspace*{2pt}
\cr
 \Psi'_{(\ltc,\ltc)}(\Gb_2,\Gb_2)\vspace*{2pt}
 \cr
 (c^{-1} \Gb_1 + \Gb_2 ,\Gb
_2)
\cr
 (\Gb_2,\Gb_2)}.
\]
Observing that $\Psi'_{(\ltc,\ltc)}$ is linear, we can conclude that
\[
c\sqrt{k} \pmatrix{
 \Psi(\ltctilde^\xi,\ltctilde) - \Psi(\ltctilde
,\ltctilde)\vspace*{2pt}
\cr
 (\ltctilde^\xi,\ltctilde) - (\ltctilde,\ltctilde)}
\weak
\pmatrix{
 \Psi'_{(\ltc,\ltc)}(\Gb_1, 0)\vspace*{2pt}
 \cr (\Gb_1,0)} =
 \pmatrix{
 \Gb_1
 \cr
  (\Gb_1,0)}.
\]
%
Continuity of the map $(\alpha, \beta, \gamma) \mapsto d(\alpha
,\beta)$ yields
\[
d \bigl( c\sqrt{k} \bigl( \Psi(\ltctilde^\xi,\ltctilde) - \Psi(\ltctilde
,\ltctilde) \bigr) ,c\sqrt{k} ( \ltctilde^\xi- \ltctilde )\bigr ) 
{\longrightarrow} 0
\]
in outer probability and thus by boundedness of the metric $d$ also in
outer expectation.
Since $c\sqrt{k} (\ltctilde^\xi- \ltctilde)\weakP\Gb_1$ we obtain
the assertion by Lemma~\ref{lem:cw}.

\subsection{\texorpdfstring{Proof of Theorem \protect\ref{theo:alphanpdm}}{Proof of Theorem 3.3}}
Again, write $T=T_i$.
We start the proof with an assertion regarding consistency of $\widehat
{\partial_p\ltc}$ and claim that for any $\delta\in(0,1)$,
%
\begin{equation}\label{consist}
\sup_{\vect{x}\in T: x_p\geq\delta} |\widehat{\partial_p\ltc
}(\vect{x})-\partial_p\ltc(\vect{x}) | \longrightarrow0
\end{equation}
in outer probability. For the proof of (\ref{consist}), split $T$ into
three subsets as indicated by its definition, and then proceed as for
the proof of Lemma 4.1 in \cite{segers2011}. The details are omitted.
Regarding the assertion of the Theorem, we set
\[
\bar\alpha^{\mathit{pdm}}_n(\vect{x}) = \beta_n(\vect{x}) - {\partial_1
\ltc}(\vect{x})\beta_{n}(x_1,\infty) - {\partial_2 \ltc}(\vect
{x})\beta_{n}(\infty,x_2).
\]
In light of Lemma~\ref{lem:cw}, it suffices to prove that $d(\alpha
_n^{\mathit{pdm}}, \bar\alpha_n^{\mathit{pdm}})$ converges to $0$ in outer probability.
By the definition of $d$, we need to show uniform convergence on the
set $T$. Because
$
| \alpha_n^{\mathit{pdm}}- \bar\alpha_n^{\mathit{pdm}}| \leq D_{n1} + D_{n2},
$
where
\[
D_{n1} = | \widehat{\partial_1\ltc} - \partial_1\ltc| |\beta
_n(\cdot,\infty)|,\qquad
D_{n2} = | \widehat{\partial_2\ltc} - \partial_2\ltc| |\beta
_n(\infty,\cdot)|,
\]
we can consider both summands $D_{np}$ separately and deal with
$D_{n1}$ as an example.
First, consider the case $\vect{x}\in[0,i]^2$. For arbitrary
$\varepsilon>0$ and $\delta\in(0,1)$,
%
\begin{equation} \label{Dnp}
\Prob^* \Bigl(\sup_{\vect{x}\in[0,i]^2} D_{n1}(\vect{x}) >\varepsilon\Bigr)
\leq
\Prob^* \Bigl(\sup_{\vect{x}\in[0,i]^2, x_1\geq\delta} D_{n1}(\vect
{x}) >\varepsilon\Bigr) +
\Prob^* \Bigl(\sup_{\vect{x}\in[0,i]^2, x_1<\delta} D_{n1}(\vect
{x})>\varepsilon\Bigr).
\end{equation}
Because $\widehat{\partial_1\ltc}$ is uniformly consistent on $\{
\vect{x}\in[0,i]^2 | x_1\geq\delta\}$, and because $\beta_n$ is
asymptotically tight
in $\ell^\infty(T)$ ($\beta_n$ converges unconditionally by the results
in Chapter 10 of \cite{kosorok2008}), the first probability on the
right-hand side
converges to 0.

Regarding the second summand, note that $F_{n1}^{-}(kx/n ) = X_{\lceil
kx \rceil:n, 1}$ (where $\lceil x\rceil=\min\{k\in\Z |\allowbreak k\geq x\}
$), so that
\[
\sup_{\vect{x}\in[0,i]^2} | \widehat{\partial_1\ltc}(\vect{x}) |
\leq\sup_{\vect{x}\in[0,i]^2, x_1\geq h} \frac{\lceil k(x_1+h)
\rceil- \lceil k(x_1-h) \rceil}{2kh} \leq1+\frac{1}{2kh} \leq2
\]
for sufficiently large $n$.
Thus the right-hand side of equation \eqref{Dnp} is bounded by
\[
\Prob^* \Bigl(\sup_{\vect{x}\in[0,i]^2, x_1<\delta} |\beta_n(\vect
{x}) | >\varepsilon/3\Bigr )
\]
eventually. Because $\beta_n\weak\Gb_{\ltctilde}$
(unconditionally), the $\limsup$ of this outer probability is bounded by
\[
\Prob\Bigl(\sup_{\vect{x}\in[0,i]^2, x_1<\delta} |\Gb_{\ltctilde
}(\vect{x})| >\varepsilon/3 \Bigr).
\]
Because $\Gb_{\ltctilde}$ has continuous trajectories and $\Gb
_{\ltctilde}(0,x_2)=0$ (almost surely),
this probability can be made arbitrary small by choosing $\delta$
sufficiently small. The case $\vect{x}\in[0,i]^2$ is completed.
For $\vect{x}\in[0,i]\times\{\infty\}$, the arguments are similar,
whereas for $\vect{x}\in\{\infty\}\times[0,i]$, we have $D_{n1}=0$,
and nothing needs to be shown. To conclude, $\sup_{\vect{x}\in T}
D_{n1}(\vect{x})$ converges to 0 in outer probability, and
because the term $\sup_{\textbf{x}\in T}D_{n_2}$ can be treated
similarly, the proof is complete.
%
%
%
%
%
%
\section{Partial derivatives of tail copulas}\label{app:B}
\begin{proposition}
The first partial derivative of a (lower or upper) tail copula $\Lambda
$ satisfies
\[
\partial_1 \Lambda(0,x) =
\cases{
\displaystyle\lim_{t\rightarrow\infty} \Lambda(1,t) & \quad if $x\in
(0,\infty)$\vspace*{2pt},
 \cr
0 &\quad if  $x =0$.}
\]
Consequently, the only tail copula that admits for continuous partial
derivatives in the origin is the tail copula corresponding to tail
independence, that is, $\Lambda\equiv0$ for either the lower or the
upper tail.
\end{proposition}

\begin{pf}
 By the groundedness and homogeneity of $\Lambda$ (see
Theorem 1 in Schmidt and Stadtm{\"{u}}ller \cite{schmstad2006}), we have
\[
\partial_1 \Lambda(0,x) = \lim_{h\rightarrow0} \frac{\Lambda
(h,x)-\Lambda(0,x)}{h} = \lim_{h\rightarrow0} \Lambda(1,x/h) = \lim
_{t\rightarrow\infty} \Lambda(1,t)
\]
for all $x\in(0,\infty)$. Similarly, $\partial_1 \Lambda(0,0) =0$.
The addendum follows by Theorem 1(iv) in Schmidt and Stadtm{\"{u}}ller \cite{schmstad2006}.
\end{pf}

As an example, note that for the Clayton copula given in
\eqref{claytoncop}, we obtain $\partial_1 \ltc(0,x) =1$ for all $\theta>0$.

\section{Auxiliary results}\label{app:C}
Here we present several technical details. We omit the proofs of the
assertions and refer the reader to B{\"{u}}cher \cite{buecher2011}, pages 103ff.
\begin{lemma}\label{BUIEunion}
Suppose that $\mathcal{G}_n$ and $\mathcal{H}_n$ are sequences of
measurable functions with envelopes $G_n$ and $H_n$, so that $(\mathcal
{G}_n,G_n)$ and $(\mathcal{H}_n,H_n)$ satisfy the bounded uniform
integral entopry condition as stated in \emph{(\ref{BUIE})}. Then the bounded
uniform entropy integral condition \emph{(\ref{BUIE})} also holds for
$\mathcal{F}_n=\mathcal{G}_n\cup\mathcal{H}_n$, with envelopes
$F_n=G_n\vee H_n$.
\end{lemma}

%

\begin{lemma}\label{bcond}
Suppose that $G_n=G_n(\vect{X}_1,\ldots\vect{X}_n,\xi_1,\ldots\xi
_n)$ is some statistic taking values in $\mathcal{B}_\infty(\bar{\R
}^2_+) $. Then a conditional version of Theorem 1.6.1 in\vspace*{-2pt} {V}an~der Vaart and Wellner \cite{vandwell1996}
holds -- namely, $G_n \weakP G$ in $\mathcal{B}_\infty(\bar{\R}^2_+)
$ is equivalent to $G_n\weakP G$ in $\ell^{\infty}(T_i)$ for every $i\in
\Nat$.
\end{lemma}
\begin{lemma}\label{lem:contboot}
Suppose that $g\dvtx\Db_1\longrightarrow\Db_2$ is a Lipschitz-continuous
map between\vspace*{-1pt} metrized topological vector spaces. If $G_n=G_n(X_1,\ldots
,X_n,\xi_1,\ldots,\xi_n)\weakP G$ in $\mathbb{D}_1$, where $G$ is\vspace*{-1pt}
tight, then $g(G_n)\weakP g(G)$ in $\Db_2$.
\end{lemma}
%
%
%
%
%
%
\begin{lemma}\label{lem:cw}
Let $Y_n=Y_n(X_1,\ldots,X_n,\xi_1,\ldots,\xi_n)$ and
$Z_n=Z_n(X_1,\ldots,X_n,\xi_1,\ldots,\xi_n)$ be\vspace*{1pt} two (bootstrap)
statistics in a metric space $(\Db,d)$,
depending on the data $X_1,\ldots,X_n$ and on some multipliers $\xi
_1,\ldots,\xi_n$. If $Y_n\weakP Y$ in $\Db$, where $Y$ is tight,
and $d(Y_n,Z_n)\Pconv0$, then also $Z_n\weakP Y$ in $\Db$.
\end{lemma}
%
%
\end{appendix}
\section*{Acknowledgements}
The authors thank Martina Stein, who typed parts of this manuscript
with considerable technical expertise. They also thank two unknown
referees and the associate editor for their constructive comments on an
earlier version of this manuscript, John Einmahl for pointing out
important references on the subject, and Johan Segers for detailed discussions.
This work was supported by the Collaborative Research Center
``Statistical modeling of nonlinear dynamic processes'' (SFB 823) of
the German Research Foundation (DFG).

\printhistory

\end{document}